\documentclass[11pt]{article}
\usepackage[english]{babel} %

\usepackage{amssymb}
\usepackage{amsmath}
\usepackage{amsthm}
\usepackage{enumitem}
\usepackage{hyperref}
\usepackage{ifthen}

\usepackage{tikz}
\usepackage[justification=centering]{caption}
\usetikzlibrary{hobby}
\usetikzlibrary{calc}
\usepackage{xfp}
\usepackage{xcolor}

\usepackage{figures}
\usepackage{figures_pent_triangl}

\advance\textheight by \topskip

\theoremstyle{plain}
\newtheorem{theorem}{Theorem}[section] 
\newtheorem*{theorem*}{Theorem} 

\newtheorem{question}[theorem]{Question}
\newtheorem{lemma}[theorem]{Lemma}
\newtheorem{corollary}[theorem]{Corollary}

\theoremstyle{definition}
\newtheorem{remark}[theorem]{Remark}

\def\Pg{\operatorname{Pg}}
\def\quote{\texttt{\char34}}
\def\Int{\operatorname{Int}}
\def\Ext{\operatorname{Ext}}

\definecolor{xblue}{RGB}{0,0,139}%
\definecolor{xgreen}{RGB}{50,205,50}%
\definecolor{xcrimson}{RGB}{220,20,60}%
\definecolor{xgold}{RGB}{255,215,0}
\definecolor{xmoccasin}{RGB}{255,228,181}

\title{Minimal pentagulations of $n$-gons}
\author{Mikhail Kabenyuk \\ \texttt{kabenyuk@gmail.com}}
\date{}

\begin{document}\maketitle

\begin{abstract}
Let $n\geq3$ be an integer.
A planar graph $G$ is called a pentagulation of an $n$-gon if all faces of $G$ are pentagons, 
except one, which is an $n$-gon.
A $3$-connected pentagulation $G$ of an $n$-gon is minimal if it has 
the smallest number of pentagons among all possible $3$-connected pentagulations of an $n$-gon.
We denote the number of pentagons of a minimal pentagulation of an $n$-gon by $\Pg(n)$.
It is known that $\Pg(3)=15$ and $\Pg(4)=14$.
We determined $\Pg(n)$ for all $3\leq n\leq12$ using computer calculations.
The calculations employed the \texttt{plantri} package, 
which generates all planar triangulations for a given number of vertices.
For each $n$ in the specified range,
we identified the number of minimal pentagulations of an $n$-gon, up to graph isomorphism.
Several open questions on this topic are presented.
The paper includes more than twenty illustrations of minimal pentagulations 
and their corresponding plane triangulations.


\end{abstract}

\section{Introduction}\label{section:Introduction}

Let a simple cycle $C_n$ of length $n$ be given in the plane.
We define a pentagulation of cycle $C_n$ (or $n$-gon) as a connected planar graph
whose faces have length five, except for the face of length $n$ bounded by cycle $C_n$.
This paper focuses on $3$-connected pentagulations.
A graph $G$ is $k$-connected if $|V(G)|>k$ and any subgraph of $G$,
obtained by removing at most $k-1$ vertices from $G$, is connected.
Note that $3$-connected planar graphs play an important role in polyhedral geometry due to Steinitz's theorem (\cite[\S51]{Steinitz}; \cite[Chapter 13]{Grunbaum}), which states that
every convex polyhedron forms a $3$-connected planar graph, and every $3$-connected planar graph can be represented as the graph of a convex polyhedron.
A convex polyhedron is defined as the intersection of a finite number of half-spaces in $\mathbb{R}^3$.

\begin{figure}[ht]
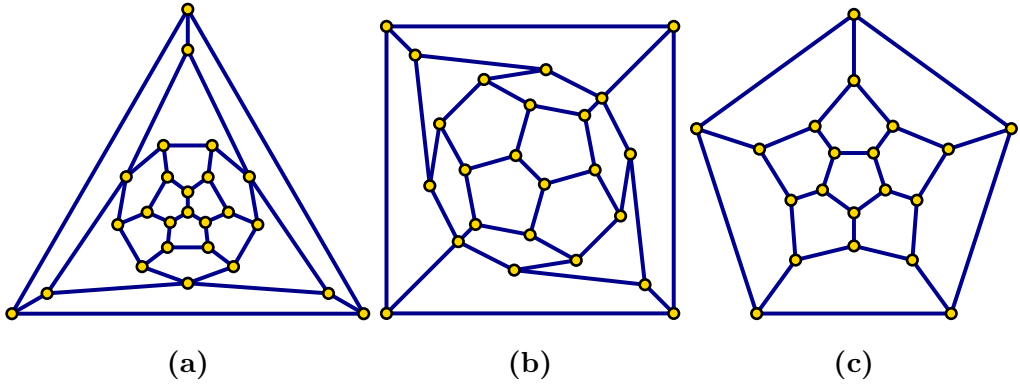

    \centering
    \TriangleMPA{27mm}{0.7mm}{(a)}
    \SquareMPOld{27mm}{0.7mm}{(b)}
    \PentagonMP{22mm}{0.7mm}{(c)}    
    \caption[Pentagulations of a triangle, a square, and a pentagon]
    {Minimal $3$-connected pentagulations of\\ a triangle ($p=15$),
    a square ($p=14$), and a pentagon ($p=11$)}
    \label{fig:TriangleSquarePentagon}
\end{figure}

Returning to pentagulations, a $3$-connected pentagulation of an $n$-gon is minimal if it contains the smallest number of pentagons among all $3$-connected pentagulations of an $n$-gon.
In answer to my question \cite{Kabenyuk} about
a minimal pentagulation of a triangle, a user known as Parcly Taxel (Jeremy Tan Jie Rui) 
constructed a $3$-connected pentagulation of a triangle with $15$ pentagons, proving that $15$ is the minimum.
Anticipating my question, he also constructed a minimal $3$-connected pentagulation of a square 
with $14$ pentagons.
Figure \ref{fig:TriangleSquarePentagon}(a,b) shows these minimal pentagulations.
In all subsequent figure captions, the number in parentheses indicates the number of pentagons.
In his computations and proofs, he used the program \texttt{plantri} \cite{Plantri}, 
specifically the module \quote \texttt{allowed\_deg.c}\quote{}
which is part of the \texttt{plantri} package.

The minimal $3$-connected pentagulation of a pentagon is the dodecahedron graph containing $11$ pentagons
(in this case, the total number of pentagons is $12$, including the outer cycle $C_5$).
In response to the question posed in \cite{licheng}
Parcly Taxel also plotted $3$-connected pentagulations of a hexagon and a heptagon
having $18$ and $21$ pentagons, respectively.
We will see below that these pentagulations are not minimal.

\begin{question}
What is the number of pentagons in the minimal $3$-connected pentagulation of an $n$-gon for $n \geq 6$?
\end{question}

We define $\operatorname{Pg}(n)$ as the number of pentagons 
in the minimal $3$-connected pentagulation of an $n$-gon. 
Thus, we know that 
$\operatorname{Pg}(3)=15$, 
$\operatorname{Pg}(4)=14$, 
$\operatorname{Pg}(5)=11$, 
while $\operatorname{Pg}(6)\leq18$, $\operatorname{Pg}(7)\leq21$.

We will prove the following statements.
\begin{theorem}\label{Th:exact_scores}
The following holds:
$$
\Pg(n)=
\left\{
  \begin{array}{ll}
    n+10, & \hbox{if $n=4,6,8$;} \\
    n+12, & \hbox{if $n=3,7,9,10,11,12$;}
  \end{array}
\right.
$$
Moreover, 
$\Pg(n)\geq n+14$ for all $n\geq13$.
\end{theorem}

\begin{theorem}\label{Th:upper_bound}
If $n\geq13$ is an integer and $n=5k+3l$ with integers $k$ and $l$ and $0\leq l\leq4$, then
\begin{equation}\label{upper_bound}
  \Pg(n)\leq2n+k+l.
\end{equation}
\end{theorem}
In addition, for $n\leq12$ we have determined the number of minimal  $3$-connected pentagulations, 
up to graph isomorphism.
\begin{theorem}\label{Th:single}
For each $n\in\{4,5,6,8,11,12\}$ there exists a unique minimal $3$-connected pentagulation of an $n$-gon.

For each $n\in\{3,7,10\}$  there are three minimal $3$-connected pentagulations of an $n$-gon.

There are four minimal $3$-connected pentagulations of a $9$-gon.
\end{theorem}

The proofs of all statements, except for Theorem \ref{Th:upper_bound}, 
rely on calculations performed with the program \texttt{plantri} \cite{Plantri}.

Unless stated otherwise, this paper considers only simple graphs,
that is, finite undirected graphs without loops or multiple edges.
We use standard notations.
If $G$ is a planar graph, then $V(G)$, $E(G)$, and $F(G)$ 
denote the sets of its vertices, edges, and faces, respectively.
Let $C$ be a simple cycle in a plane graph. 
Then the rest of the plane is
partitioned into two disjoint open sets called the interior and exterior of $C$.
We shall denote the interior and exterior of $C$ by $\Int(C)$ and $\Ext(C)$, respectively,
and we do not use any special notation for their closures
(however, see \cite[Chapter 9]{Bondy} for other conventions).
Throughout this paper, references to the number of graphs (triangulations, pentagulations, at al) 
with a given property, we refer to the number of distinct isomorphism classes.

The paper is organized as follows.
Section \ref{section:Preliminaries} provides preliminary information,
including the introduction of the parameter $\Delta^+(G)$ for a pentagulation $G$.
In Sections \ref{section:Special_triangulations} and 
\ref{section:Triangulations with dfc(G)=1} 
we prove that there are no plane pentagulations $G$ for which $\Delta^+(G)=0$ or $\Delta^+(G)=1$.
Section \ref{section:Triangulations with dfc(G)=2} presents constructions of pentagulations $G$ with $\Delta^+(G)=2$, 
while 
Section \ref{section:Triangulations with dfc(G)=3} examines pentagulations for which  $\Delta^+(G)=3$.
Section \ref{section:Triangulations with dfc(G)=3} concludes by proving Theorems \ref{Th:exact_scores} and \ref{Th:single}.
In Section \ref{section:Upper bound}, we prove Theorem \ref{Th:upper_bound}.
Finally, Section \ref{section:Questions} formulates several open questions.
The appendix (Section \ref{section:Appendix}) includes a collection of diagrams of $3$-connected pentagulations $G$
with $\Delta^+(G)\geq3$.

\section{Preliminaries}\label{section:Preliminaries}
\begin{figure}[ht]
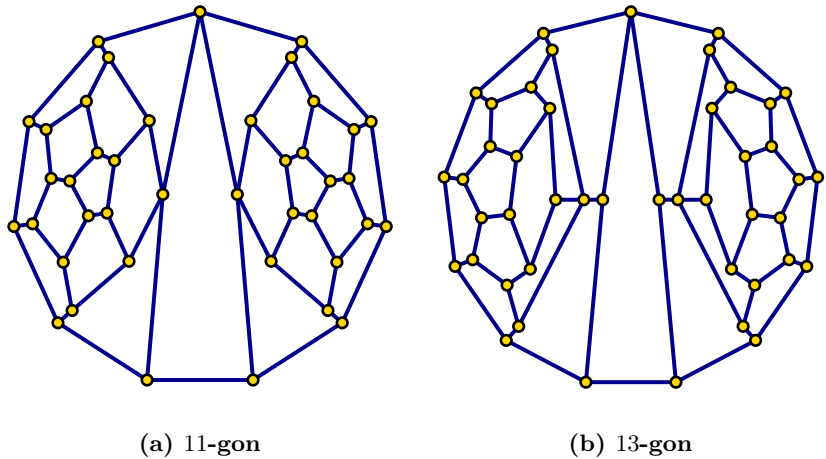

    \centering
    \UndecagonS{25mm}{0.7mm}{\footnotesize (a) $11$-gon}
    \quad    
    \TridecagonS{25mm}{0.7mm}{\footnotesize (b) $13$-gon}
    \caption[2-connected pentagulations of an undecagon and a tridecagon]
    {2-connected pentagulations of an undecagon and a tridecagon}
    \label{fig:UndecagonTridecagon_2-connected}
\end{figure}

Note that if the pentagulation is $3$-connected, then every vertex has a degree of at least three.
However, the converse is not true in general 
For instance, consider the pentagulations of an undecagon (a) and a tridecagon (b) are shown in Figure \ref{fig:UndecagonTridecagon_2-connected}.
Nevertheless, if every vertex in a $2$-connected planar graph has a degree of at least three,
and the intersection of any two face boundaries is connected, then this graph is $3$-connected
\cite[Theorem 1.14 and Proposition 1.10]{Dunlaing}.
It is worth noting that this statement will not be used in the paper.

Let $G$ be a $3$-connected planar graph with $|F(G)|=p+1$ faces, 
$p$ of which are pentagons and one an $n$-gon
(i.e., $G$ is a pentagulation of an $n$-gon), where 
$n\geq3$ and $n\neq5$. 
The number of edges in $G$ is given by
\begin{equation}\label{formula:NumberEdgesGraph}
    |E(G)|=\frac{5p+n}{2}.
\end{equation}
It follows that the integers $p$ and $n$ must have the same parity.
By Euler's formula for planar graphs $|V(G)|-|E(G)|+|F(G)|=2$, and therefore
\begin{equation}\label{formula:NumberVerticesGraph}
    |V(G)|=\frac{3p+n+2}{2}.
\end{equation}
Since every vertex of $G$ has a degree of at least three, 
it follows from the handshake lemma that
\begin{equation}\label{formula:DeficiencyDegreesGraph}
    \Delta^+(G)=\sum_{v\in V(G)}(\deg(v)-3)=2|E(G)|-3|V(G)|=\frac{p-n}{2}-3\geq0.
\end{equation}
We will call the number $\Delta^+(G)$ \textit{the degree surplus} of the graph $G$.

Now let $G^*$ be the dual of $G$.
As we know, the dual of a $3$-connected planar graph is 
also a $3$-connected planar graph \cite[Theorem 10]{Whitney}.
Furthermore, by Whitney's theorem (\cite[Theorem 11]{Whitney}, \cite[Theorem 4.3.2]{Diestel}),
the dual is uniquely determined and does not depend on the embedding of $G$ in the plane.
In addition, the following holds:
\begin{equation}\label{formula:NumberVerticesDual}
    |V(G^*)|=|F(G)|=p+1,
\end{equation}
\begin{equation}\label{formula:NumberFacesDual}
    |F(G^*)|=|V(G)|=\frac{3p+n+2}{2},
\end{equation}
\begin{equation}\label{formula:NumberEdgesDual}
    |E(G^*)|=|E(G)|=\frac{5p+n}{2}.
\end{equation}
The graph $G^*$ has $p$ vertices of degree $5$ and exactly one vertex of degree $n$.
If $n=5$, all the formulas above remain valid.

If $G$ is a $3$-connected pentagulation of an $n$-gon with
$\Delta^+(G)=0$, then $G^*$ is a plane triangulation; 
that is, every face of $G^*$, including the outer face, is bounded by a triangle and $p=n+6$.
An example of such a triangulation is the icosahedral graph 
(Figure \ref{fig:TriangulationTriangleSquarePentagon}(a)) 
whose dual is the dodecahedral graph
(Figure \ref{fig:TriangleSquarePentagon}(c)).
The equality $\Delta^+(G)=0$
implies that the dodecahedron is a minimal $3$-connected pentagulation of a pentagon.

In what follows, we focus on graphs that are dual to pentagulations.
For example, the dual $G^*$ of the pentagulation of a triangle 
(Figure \ref{fig:TriangleSquarePentagon}(a)) has $16$ vertices: 
$15$ of degree $5$ and one of degree $3$.
Three of the faces in $G^*$ are squares,  while the remaining ones are triangles.
If we prove that no $3$-connected planar graphs exist with fewer than $16$ vertices,
where all but one have degree $5$ and one has degree $3$,
this would imply that the pentagulation of the triangle shown
in Figure \ref{fig:TriangleSquarePentagon}(a) is indeed minimal.

A triangulation of a cycle $C_n$ (or an $n$-gon) is a connected planar graph
whose faces are all triangular, except for one, which is bounded by a cycle of length $n$.
Hereafter, we assume that the face bounded by the cycle $C_n$ is the external face, and
all other triangular faces are in the interior of $C_n$.
A separating $n$-cycle in a plane graph is an $n$-cycle such that both the
interior and the exterior contain at least one vertex.
A $3$-cycle in a plane graph is called an innermost separating $3$-cycle if its interior
contains no other separating $3$-cycles.
It is clear that an innermost separating $3$-cycle exists if
there is at least one separating $3$-cycle.

In the following, we derive a formula for the number of interior points of a cycle triangulation.
\begin{lemma}
\label{lemma:on interior points}
  Let $G$ be a plane triangulation and $v\in V(G)$.
  All neighbours of $v$ form a cycle $C$ of length $n$. 
  We assume that the vertex $v$ lies in $\Ext(C)$, and 
  all other vertices lie in $C$ or in $\Int(C)$.
  Let $A$ be the set of internal vertices of $C$. 
  Let $I(v)=|A|$.  
  
  We further denote $p_v$ as the sum of the degree excesses of all vertices other than 
  $v$ over $5$, i.e. 
  $$
  p_v=\sum_{x\in V(G),\ x\neq v}\deg(x).
  $$
  Then  
  \begin{equation}
  \label{eq:number of interior points}
        I(v)=p_v+6.
  \end{equation}
\end{lemma}
\textit{Proof.}
We have 
$$
|V(G)|=n+I(v)+1,\ 2|E(G)|=3|F(G)|. 
$$
Using Euler’s formula for planar graphs, we obtain from these equations
$$
2=n+I(v)+1-|E(G)|+F(G)=n+I(v)+1-|E(G)|/3.
$$
Hence, 
\begin{equation}
    \label{eq:I(v)=|E(G)|/3-n+1}
    I(v)=|E(G)|/3-n+1.
\end{equation}
By the handshaking lemma, 
\begin{align*}
  2|E(G)|&=\sum_{x\in V(G)}\deg(x) \\
  &=\sum_{x\in V(G),\ x\neq v}\deg(x)-5(|V(G)|-1)+5(|V(G)|-1)+n\\
  &=p_v+6n+5I(v).
\end{align*}
Substituting the expression for $2|E(G)|=p_v+6n+5I(v)$ into the formula (\ref{eq:I(v)=|E(G)|/3-n+1}), 
we obtain the formula (\ref{eq:number of interior points}).
The lemma is proved.

We finish this section by proving a statement about the chords of an $n$-gon in a plane graph.
\begin{lemma}[Existence of chords]
\label{lemma:Existence of chords}
  Let $n\geq3$.
  Let $G$ be a plane graph and let $C$ be an $n$-gon face of $G$.
  If $\Int(C)$ contains no vertices of $G$,
  then there exist $n-3$ chords of $C$ lying in $\Int(C)$,
  having no common points within $\Int(C)$, and subdividing $\Int(C)$ into $n-2$ triangles.
\end{lemma}
\textit{Proof.}
Recall that a chord of the $n$-gon $C$ is a Jordan curve
connecting two non-consecutive  vertices of $C$ and lying entirely within $\Int(C)$.

We prove the lemma by induction on $n$.
It is true when $n=3$.
Suppose that it holds when $n=N-1\geq3$,
and let $C$ be an $N$-gon with vertices $x_1,x_2,\ldots,x_N$.

If $x_1$ and $x_3$ are not adjacent vertices in the graph $G$,
then we can draw the chord $L$ connecting the vertices $x_1$ and $x_3$ and lying within $\Int(C)$.
By the induction hypothesis, the $(N-1)$-gon $C'=x_1x_3x_4\ldots x_N$ has $N-4$ chords
that meet the necessary requirements.
Since every chord of $C'$ is also a chord of $C$,
it follows that, together with $L$, we have $N-3$ chords of the $N$-gon $C$.

\begin{figure}[ht]
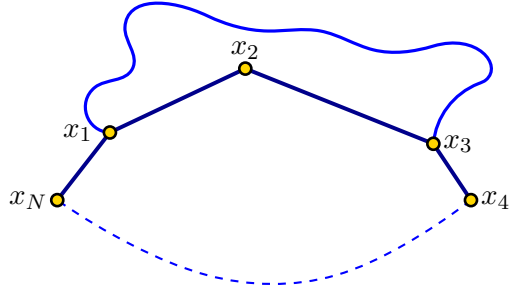

    \centering
    \VerticesAdjacentInG
    \caption[Chords of an n-gon]
    {The vertices $x_1$ and $x_3$ are adjacent in the graph $G$}
    \label{fig:Chords of an n-gon}
\end{figure}

Now consider the case where $x_1$ and $x_3$ are adjacent in the graph $G$ (see Figure \ref{fig:Chords of an n-gon}).
By the Jordan curve theorem, the vertex $x_2$ cannot be adjacent to any
vertex $x_i$ with $i=4,\ldots,N$ in $G$.
Consequently, each of the vertices $x_4,\ldots,x_N$ can be connected to $x_2$ by chords.
Thus, the result holds for $n=N$, and by the principle of induction,
it is valid for all $n\geq3$. 
Thus, the lemma is proved.

\section{Triangulations with all but one vertex of degree 5}
\label{section:Special_triangulations}
This section is devoted to proving the following statement.
\begin{theorem}\label{Th:pure_triangulation}
    Let $n\geq3$ be an integer with $n\neq5$.
    There does not exist a $3$-connected plane triangulation 
    in which exactly one vertex has degree $n$, while all other vertices have degree $5$.
\end{theorem}

\begin{remark}
\label{remark:About a theorem Grunbaum}
  When $n$ is not a multiple of $5$, Theorems 
  \ref{Th:pure_triangulation} and 
  \ref{Th:No triangulations with dfc=1}
  in the next section follow as 
  corollaries of Theorem 13.4.2(5) in \cite{Grunbaum}.
  However, Gr$\rm\ddot{u}$nbaum's theorem does not cover  
  the case where $n$ is a multiple of $5$ and $n\neq5$.  
\end{remark}

\begin{corollary}\label{corollary:LowerBound:n+8}
    Let $n\geq3$ be an integer and $n\neq5$. 
    Any $3$-connected pentagulation of an $n$-gon contains at least $n+8$ pentagons.
\end{corollary}
\textit{Proof.}
Let $p$ be the number of pentagons of a $3$-connected pentagulation.
From (\ref{formula:DeficiencyDegreesGraph}), $p\geq n+6$.
By Theorem \ref{Th:pure_triangulation}, $p\neq n+6$.
Thus, $p\geq n+8$, as $p$ and $n$ have the same parity.
The corollary is proved.

Before proving Theorem \ref{Th:pure_triangulation}, we establish a few auxiliary results.
\begin{lemma}\label{lemma:Triangulation_n-cycle}
    Let $H$ be a plane triangulation of an $n$-cycle $C$.
    Let $A$ be the set of vertices inside $C$, where each vertex in $A$ has a degree at least $5$.
    \begin{enumerate}[label=$\arabic*.$,align=left,left=\the\parindent]
    \item
    If $k=|A|$, then
    \begin{equation}\label{number of edges of the triangulation cycle}
        |E(H)|=3k+2n-3.
    \end{equation}
    \item
    If $s$ is the sum of the degrees of the vertices in $C$, then
    \begin{equation}\label{number of internal vertices in the triangulation of the n-cycle}
        k\geq s-4n+6.
    \end{equation}
    \item
    If $n=3$ or $n=4$ with $k\geq1$, or $n=5$ with $k\geq2$, then
    \begin{equation}\label{weak boundary for the innermost}
        k\geq9-n.
    \end{equation}
    \item
    If in addition the subgraph induced by the set $A$ is not a triangulation of the plane, then
    \begin{equation}\label{strict boundary for the innermost}
        k\geq10-n.
    \end{equation}
    \end{enumerate}
\end{lemma}
\textit{Proof.}
1. Since
$|V(H)|=k+n$ and $2|E(H)|=3(|F(H)|-1)+n$,
it follows from Euler's theorem for planar graphs that
(\ref{number of edges of the triangulation cycle}) holds.
Note that this result does not require any conditions on the degrees of the vertices in $A$.

2. Since $2|E(H)|\geq 5k+s$, it follows from (\ref{number of edges of the triangulation cycle}) that
the inequality (\ref{number of internal vertices in the triangulation of the n-cycle}) holds.

3. Let $H'$ be a graph induced by the set $A$ and $t=|E(H')|$.
If $v\in A$, then
the number of edges connecting $v$ to vertices in $C$ is at least $5-\deg_{H'}(v)$.
Therefore, the number of edges $r$ joining vertices of $A$ with vertices of $C$ satisfies the inequality
$$
r\geq\sum_{v\in A}5-\deg_{H'}(v)=5k-2t.
$$
Clearly, $s\geq2n+r$ (with equality if no edges of $H$ are chords of $C$), and thus 
\begin{equation}\label{boundary for s}
s\geq2n+5k-2t.
\end{equation}
It is easy to verify that if $n=3,4$ and $k\geq1$ or $n=5$ and $k\geq2$, then $k=|A|\geq3$.
Hence for the planar graph $H'$ the inequality $t\leq3k-6$ holds,
so it follows from (\ref{boundary for s}) that
\begin{equation}\label{boundary for s at k>2}
s\geq2n-k+12.
\end{equation}
Combining inequalities (\ref{number of internal vertices in the triangulation of the n-cycle})
and (\ref{boundary for s at k>2}), we derive
$$
k\geq-k-2n+18\ \Rightarrow\ k\geq9-n,
$$
and the inequality (\ref{weak boundary for the innermost}) is proved.

4. To prove the inequality (\ref{strict boundary for the innermost}), 
it is sufficient to note that if the planar graph $H'$  has 
at least three vertices and is not a triangulation of the plane,
then $|E(H')|<3|V(H')|-6$ or $t<3k-6$.
The lemma is proved.

\begin{figure}[ht]
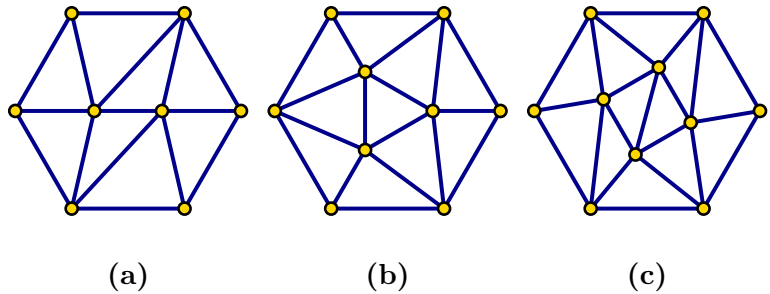

    \centering
    \TriangulationsHexagon
    \caption[Triangulations of a hexagon]
    {Triangulations of a hexagon containing interior two, three, \\and four vertices of degree $5$}
    \label{fig:TriangulationsHexagon}
\end{figure}

\begin{remark}
  As shown in Figure \ref{fig:TriangulationsHexagon}(a), 
  the inequality $(\ref{weak boundary for the innermost})$ does not hold for $n=6$.
\end{remark}

\begin{corollary}\label{corollary:interior_of_a_n-cycle (weak variant)}
    Let $H$ be a triangulation of an $n$-cycle $C$.
    Let $A$ be the set of vertices of $H$ lying in the interior of $C$, and let
    $k=|A|$.
    Assume every vertex in $A$ has a degree of at least $5$. 
    Then
    \begin{enumerate}[label=$\arabic*.$,align=left,left=\the\parindent]
      \item if $n=3$ and $k>0$, then $k\geq7$;
      \item if $n=4$ and $k>0$, then $k\geq6$;
      \item if $n=5$ and $k>1$, then $k\geq5$.
    \end{enumerate}
\end{corollary}
\textit{Proof.}
Let $H'$ be the same as in the proof of Lemma \ref{lemma:Triangulation_n-cycle}.
If $H'$ is not a triangulation of the plane, 
then all inequalities follow directly from (\ref{strict boundary for the innermost}).
Let $H'$ be a triangulation of the plane.
Let $T$ be a triangle in $H'$ containing all vertices of $A$ and 
and let $A'$ denote the set of interior vertices of $T$.
From (\ref{weak boundary for the innermost}), we deduce that $|A|>3$, 
and hence $A'\neq\varnothing$.
Applying inequality (\ref{weak boundary for the innermost}) to the $3$-cycle $T$, we get
$|A|=|A'|+3\geq9$.
Stronger inequalities hold in this case, and these will be needed later.
The corollary is proved.

\begin{figure}[ht]
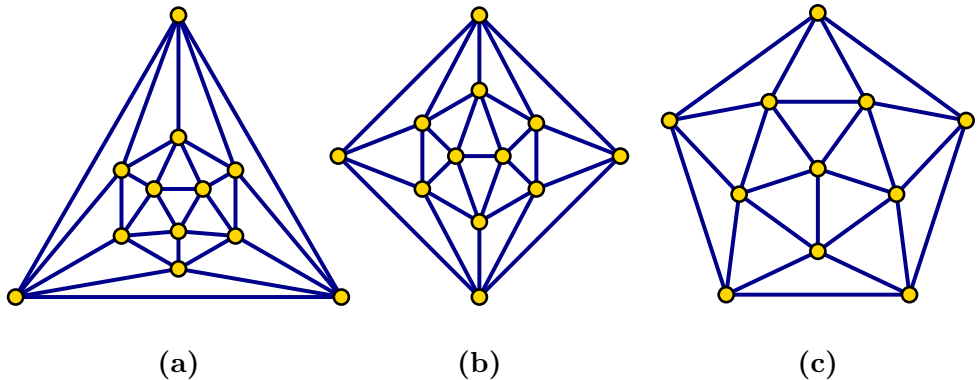

    \centering
    \TriangulationTriangleSquarePentagon{25mm}{1mm}
    \caption[Triangulations of a triangle (Icosahedral graph), a square, and a pentagon]
    {Triangulations of a triangle, a square, and a pentagon}
    \label{fig:TriangulationTriangleSquarePentagon}
\end{figure}

A stronger version of Corollary \ref{corollary:interior_of_a_n-cycle (weak variant)} is given in the next lemma.
The bounds given in this lemma are attained by the icosahedral graph 
(Figure \ref{fig:TriangulationTriangleSquarePentagon}(a)), as well as by the graphs in Figures \ref{fig:TriangulationTriangleSquarePentagon}(b) and \ref{fig:TriangulationTriangleSquarePentagon}(c).

\begin{lemma}\label{lemma:interior_of_a_n-cycle (reinforced version)}
    Let $H$ be a triangulation of an $n$-cycle $C$.
     Let $A$ denote the set of vertices of $H$ lying in the interior of $C$, 
    and let $k=|A|$.
    Assume every vertex in $A$ has a degree of at least $5$. Then
    \begin{enumerate}[label=$\arabic*.$,align=left,left=\the\parindent]
      \item if $n=3$ and $k>0$, then $k\geq9$;
      \item if $n=4$ and $k>0$, then $k\geq8$;
      \item if $n=5$ and $k>1$, then $k\geq6$;
      \item There exists exactly one triangulation of a pentagon with six vertices in its interior.
    \end{enumerate}
\end{lemma}
\textit{Proof.}
According to Corollary \ref{corollary:interior_of_a_n-cycle (weak variant)},
it suffices to show that the following cases are impossible:
$n=3$, $k=7,8$; $n=4, k=6,7$; $n=5, k=5$.
We will consider these cases in reverse order.

\begin{figure}[ht]
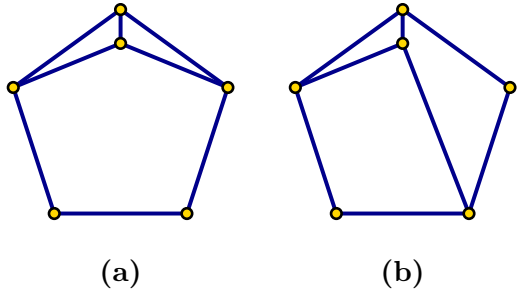

    \centering
    \PartitioningCycleOrderFive
    \caption[Partitioning the interior of a $5$-cycle]
    {Partitioning the interior of a $5$-cycle}
    \label{fig:PartitioningCycleOrderFive}
\end{figure}

\begin{figure}[ht]
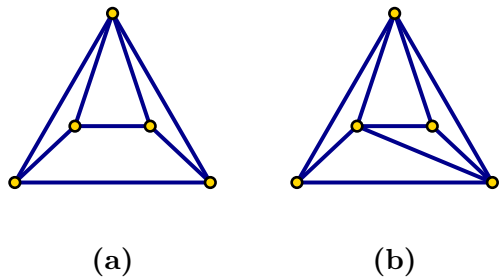

    \centering
    \GraphsOrderFive
    \caption[Two graphs of order $5$]
    {Two graphs of order $5$}
    \label{fig:Two graphs of order 5}
\end{figure}

\textsc{Case } $n=5$, $k=5$.
Let $H'$ be the same as in the proof of Lemma \ref{lemma:Triangulation_n-cycle}.
If $v\in A$ and $\deg_{H'}(v)\leq2$, then
the number of edges connecting $v$ with the vertices from $C$
is at least three.
For instance, if three edges connect $v$ to the vertices of $C$, 
then either a triangle, a quadrilateral, or a pentagon 
must contain two to four interior vertices 
Then, by Corollary \ref{corollary:interior_of_a_n-cycle (weak variant)} 
$k-1\geq6$, which contradicts the assumption that $k=5$.
(see Figure \ref{fig:PartitioningCycleOrderFive}(a) and \ref{fig:PartitioningCycleOrderFive}(b)).
Note that the same reasoning implies that the chords of $C$ cannot be edges of $H$. 
We will need this just below.
So $\deg_{H'}(v)\geq3$ for every $v\in A$.
Since the graph $H'$ has $5$ vertices,
it follows that $H'$ has at least $8$ edges.
Moreover, since $H'$ is a planar graph,
there are exactly two possibilities, as shown in Figure
\ref{fig:Two graphs of order 5}(a) and \ref{fig:Two graphs of order 5}(b).
They are obtained from the complete graph $K_5$ by the removal of one or two non-adjacent edges.
Both graphs are $3$-connected, 
so they actually have exactly one embedding in the plane (up to isomorphism).
For a plane representation, this means that either a triangle or a quadrilateral has one or two interior vertices, which contradicts Corollary \ref{corollary:interior_of_a_n-cycle (weak variant)}.

\begin{figure}[ht]
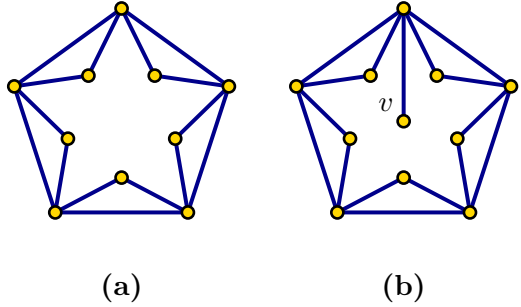

    \centering
    \UniquenessTriangulationPentagon
    \caption[Uniqueness of the triangulation of a pentagon]
    {Uniqueness of the triangulation of a pentagon \\ with six interior vertices}
    \label{fig:Uniqueness of the triangulation of a pentagon}
\end{figure}

We will prove that there is exactly one triangulation of a pentagon with six interior vertices.
This triangulation is shown in Figure \ref{fig:TriangulationTriangleSquarePentagon}(c).
We know that $A$ has no vertices adjacent to three or more vertices of $C$,
that the chords of $C$ are not edges of $H$, and 
that every edge of $C$ belongs to a triangular face in $H$.
Consequently, each vertex of $C$ is incident to at least two edges connecting it to the interior of $C$.
Figure \ref{fig:Uniqueness of the triangulation of a pentagon}(a) illustrates this situation.
If the sixth vertex $v\in A$ is adjacent to a vertex of $C$ 
(see Figure \ref{fig:Uniqueness of the triangulation of a pentagon}(b)), 
then $s\geq21$, where $s$ denotes the sum of the degrees of the vertices in $C$
From (\ref{number of internal vertices in the triangulation of the n-cycle}) we get $k\geq7$. 
This This contradicts the assumption that $k=5$.
Thus, $v$ is adjacent only to the five interior vertices, resulting in the graph shown in Figure \ref{fig:TriangulationTriangleSquarePentagon}(c).
Note that all vertices of the single pentagonal face in this graph have degree $4$.

\textsc{Case} $n=4$, $k=6,7$.
Consider a $4$-cycle $x_1x_2x_3x_4$. 
Choose $v\in A$ so that no vertices of $A$ lie in the interior of the triangle $x_1vx_2$.
By the previous case, the interior of the $5$-cycle $vx_2x_3x_4x_1$ contains at least six vertices.
If there are more than $6$ vertices, then $k\geq8$, which contradicts the assumption.
Suppose that the interior of this $5$-cycle contains exactly six vertices.
Also according to the previous item,
there exists exactly one triangulation of a $5$-cycle  $vx_2x_3x_4x_1$ 
with six interior vertices of degree $5$ or higher.
At that, the vertices of this $5$-cycle have degree $4$.
Therefore $\deg(v)=4$. This is a contradiction.

\textsc{Case} $n=3$, $k=7,8$.
The proof of this case relies on the proof of the previous case.
Consider a $3$-cycle $x_1x_2x_3$, and 
choose $v\in A$ so that no vertices of $A$ lie inside the triangle $x_1vx_2$.
The interior of the $4$-cycle $vx_2x_3x_1$ contains at least eight vertices; hence $k\geq9$.
So the lemma is completely proved.

\begin{figure}[ht]
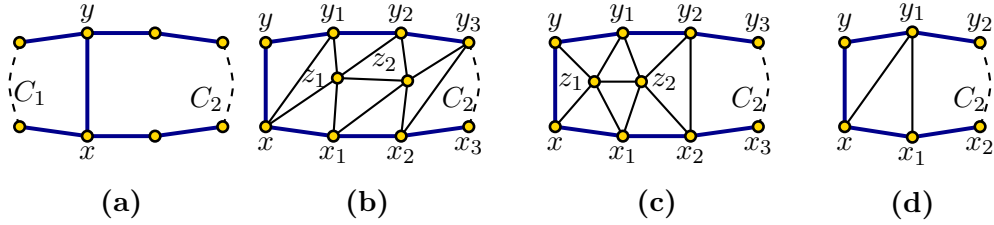

    \centering
    \NoChordsLemma{0.9}{3.5}{0.07}
    \caption[Illustration to the no chords lemma]
    {Illustration to the no chords lemma}
    \label{fig:NoChordsLemma}
\end{figure}

\begin{lemma}[The no chords lemma]
\label{lemma:The no chords}
    Let $H$ be a triangulation of an $n$-cycle $C$ and $V(H)\subset\Int(C)\cup V(C)$.
    Assume that one of the following two conditions is satisfied:
    \begin{enumerate}[label=$(\roman*)$,align=left,left=\the\parindent]
      \item each vertex in $C$ has degree $4$, except for at most one vertex of degree $5$;
      \item each vertex in $C$ has degree $4$, except for exactly two adjacent vertices of degree $5$, and each vertex in $\Int(C)$ has degree $5$.
    \end{enumerate}
    Then no chord of the cycle $C$ is an edge of the graph $H$.
\end{lemma}
\textit{Proof.}
Suppose, on the contrary, that there is
an edge in $H$ connecting two vertices $x,y\in V(C)$ that are not consecutive in $C$;
that is, $xy$ is a chord of the cycle $C$ 
(see Figure \ref{fig:NoChordsLemma}(a)).
The edge $xy$ divides $C$ into two cycles $C_1$ and $C_2$, each of length at least three.
Let $H_i$ be the subgraph of $H$ induced by the set 
$(V(C_i)\cup\Int(C_i))\cap V(H)$ for $i=1,2$.

Without loss of generality, we can assume that if the vertices of degree $5$ are the ends of a chord of $C$, 
then they are $x$ and $y$.
Given this assumption, one of the cycles, say $C_2$, does not contain vertices of degree $5$ other than $x$ and $y$.

We denote by $\deg_i(v)$ the degree $v\in V(H_i)$, and $\deg(v)=\deg_H(v)$ is the degree in $H$.
Clearly, $\deg_i(v)=\deg(v)$ for any $v\in V(H_i)$ if $v\neq x$ and $v\neq y$. 
Therefore, we use the notation $\deg_i$ only for $x$ and $y$.

If for some $i\in\{1,2\}$ we have $\deg_i(x)=\deg_i(y)=2$, 
then the cycle $C_i$ is a triangle $xyz$, where $z\in V(C_i)$ and $H_i=C_i$. 
Hence $\deg(z)=2$, which contradicts the condition of the lemma. 
Consequently, the following inequality holds for each $i=1,2$:
\begin{equation}\label{ineq:deg_i(x)+deg_i(y)>=5}
  \deg_i(x)+\deg_i(y)\geq5.
\end{equation}
On the other hand, for each $i=1,2$, we have:
\begin{equation}\label{ineq:deg_i(x)+deg_i(y)<=7}
    \deg_i(x)+\deg_i(y)\leq7,
\end{equation}
since $\deg(x)\leq5$, $\deg(y)\leq5$, which implies $\deg_i(x)\leq4$, $\deg_i(y)\leq4$, 
and in the left part of inequality
(\ref{ineq:deg_i(x)+deg_i(y)<=7}) the edge $xy$ is counted twice. 
Moreover, if, for instance, $\deg_2(x)+\deg_2(y)=7$, and 
$\deg_2(x)=4$, $\deg_2(y)=3$, then $\deg(x)=5$ and $\deg(y)=4$ or $5$.
If $\deg(y)=4$, then $\deg_1(x)=\deg_1(y)=2$ which, as we have seen, is impossible.
If $\deg(y)=5$, then $\deg_1(x)+\deg_1(y)=5$.   
Since in this case both cycles $C_1$ and $C_2$ have no vertices of degree $5$ other than $x$ and $y$, 
we can swap $C_1$ and $C_2$.
Given this and inequality (\ref{ineq:deg_i(x)+deg_i(y)>=5}), 
we can assume that the following inequalities are satisfied
\begin{equation}\label{ineq:5<=deg_2(x)+deg_2(y)<=6}
  5\leq\deg_2(x)+\deg_2(y)\leq6.
\end{equation}
From now on, we assume that inequalities (\ref{ineq:5<=deg_2(x)+deg_2(y)<=6}) hold 
and all vertices of the cycle $C_2$ have degree $4$ except for $x$ and $y$.
Since we can always interchange $x$ and $y$, there are only the following three cases due to (\ref{ineq:5<=deg_2(x)+deg_2(y)<=6}):
\begin{enumerate}[label=$\arabic*.$,align=left,left=\the\parindent]
  \item $\deg_2(x)=4$, $\deg_2(y)=2$;
  \item $\deg_2(x)=3$, $\deg_2(y)=3$;
  \item $\deg_2(x)=3$, $\deg_2(y)=2$.
\end{enumerate}

\textsc{Case 1:} $\deg_2(x)=4$, $\deg_2(y)=2$.

Since $\deg_2(y)=2$,
the third vertex $y_1$ of a triangular face incident with $xy$ lies on $C_2$,
and $y$ and $y_1$ are consecutive vertices on $C_2$ 
(see Figure \ref{fig:NoChordsLemma}(b)).

Because $\deg_2(x)=4$ and $\deg(y_1)=4$, the triangle incident with $xy_1$ (different from $xyy_1$)
must have its third vertex lying in $\Int(C_2)$; 
we denote this vertex by $z_1$.
For the same reason, the triangles incident with edges $xz_1$ and $y_1z_1$ 
and distinct from the previously considered ones, 
are $xz_1x_1$ and $y_1z_1y_2$, respectively, 
where $x$, $x_1$ and $y_1$, $y_2$ are consecutive vertices on $C_2$.

Now, since $\deg(z_1)=5$, $\deg(x_1)=4$, and $\deg(y_2)=4$, 
the new triangles incident with edges $x_1z_1$ and $z_1y_2$ are $x_1z_1z_2$ and $y_2z_1z_2$, respectively.
where $z_2\in\Int(C_2)$.

Proceeding similarly, since $\deg(x_1)=4$ and $\deg(y_2)=4$ 
the other triangles incident with edges $x_1z_2$ and $y_2z_2$ are
$x_1z_2x_2$ and $y_2z_2y_3$, respectively, 
with $x_1$ and $x_2$, and $y_2$ and $y_3$ being consecutive vertices on $C_2$.

Continuing this pattern, and since $\deg(z_2)=5$ we find that $x_2y_3$ must be an edge in $H$.
This leads to the existence of edge $x_3y_3$, where $x_3$ and $x_2$ are consecutive vertices on $C_2$.
Consequently, $\deg(x_3)=2$, which contradicts the condition that all vertices on $C_2$ 
(except possibly $x$ and $y$) have degree $4$.

\textsc{Case 2:} $\deg_2(x)=3$, $\deg_2(y)=3$.

In this case, the third vertex $z_1$ of the face incident with $xy$ lies in $\Int(C_2)$
(see Figure \ref{fig:NoChordsLemma}(c)).
Then, the triangles incident with edges $xz_1$ and $yz_1$, different from $xyz_1$, 
are $xz_1x_1$ and $yz_1y_1$, respectively,
where $x$ and $x_1$, and $y$ and $y_1$ are consecutive vertices of $C_2$.
Continuing as in \textsc{Case 1}, we eventually reach $x_3=y_3$ and $\deg(x_3)=2$, 
which contradicts the condition.

\textsc{Case 3:} $\deg_2(x)=3$, $\deg_2(y)=2$.

This case is simpler than the previous ones
(see Figure \ref{fig:NoChordsLemma}(d)).
Since $\deg_2(y)=2$,
the triangle incident with edge $xy$ is $xyy_1$, 
where $y_1$ is the vertex adjacent to $y$ on $C_2$.
Because $\deg_2(x)=3$,
the triangle incident with edge $xy_1$ and different from $xyy_1$ is $xx_1y_1$
where $x_1$ is the vertex adjacent to $x$ on $C_2$.

Next, since $\deg(x_1)=\deg(y_1)=4$, 
we consistently obtain that $H$ has the edges $x_1y_2,x_2y_2$,
where $x_2$ and $y_2$ are the vertex of $C_2$ adjacent to $x_1$ and $y_1$, respectively.

Continuing this, we consistently obtain sequences of vertices
$x,x_1,x_2,\ldots$ and $y,y_1,y_2,\ldots$ along $C_2$.

Since $C_2$ is a finite cycle, this process cannot continue indefinitely. 
Obviously, as soon as this process breaks down, 
we arrive at a contradiction 
that the degree of every vertex on $C_2$, except $x$ and $y$, is $4$.

In either case, we arrive at a contradiction. 
Therefore, no chord of the cycle $C_2$ is an edge of the graph $H$,
the lemma is proved.

\begin{remark}
    \label{remark:case 3 of Lemma no chords}
    In Case 3 of Lemma \ref{lemma:The no chords} 
    there are no specific conditions imposed on the vertices lying in $\Int(C_2)$. 
    Therefore, the following statement holds:

    Let $n\geq3$ be an integer.   
    There is no triangulation of an $n$-cycle $C$ such that 
    all vertices of $C$ have degree $4$,
    except for two consecutive vertices $x$ and $y$ on $C$, for which $\deg(x)+\deg(y)\leq5$. 
        
    Indeed, 
    if $\deg(x)=2$ and $\deg(y)=2$, 
    then $C$ must be a triangle formed by the vertices $x$, $y$, and a third vertex $z\in V(C)$. 
    In this case, $\deg(z)=2$, which contradicts the assumption that all other vertices have degree $4$.
      
    If $\deg(x)=2$ and $\deg(y)=3$ or $\deg(x)=3$ and $\deg(y)=2$, 
    then this situation corresponds exactly to Case 3 of Lemma \ref{lemma:The no chords}, 
    which states that such a triangulation does not exist.
\end{remark} 

\begin{figure}[ht]
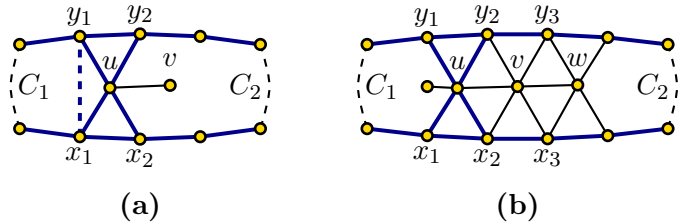

    \centering
    \CommonVertexTriangles{0.8}{3.5}{0.07}
    \caption[Common vertex of triangles]
    {Common vertex of triangles}
    \label{fig:CommonVertexTriangles}
\end{figure}

\begin{lemma}\label{lemma:no common vertex of triangles}
    Let $H$ be a triangulation of an $n$-cycle $C$, where $V(H)\subset\Int(C)\cup V(C)$.    
    Let us assume that one of the following three conditions is satisfied:
    \begin{enumerate}[label=$(\roman*)$,align=left,left=\the\parindent]
      \item       
      Each vertex in $\Int(C)$ has degree $5$, 
      and at most one or two adjacent vertices in $C$ have degree $5$,      
      while the rest of the vertices in $C$ have degree $4$.
      
      \item 
      Each vertex in $\Int(C)$ has degree $5$, except for at most one vertex of degree $6$, and
      each vertex in $C$ has degree $4$, except for at most one vertex of degree $5$.
      
      \item 
      At most two vertices in $\Int(C)$ have degree $6$, while all other vertices in $\Int(C)$ have degree $5$, and each vertex in $C$ has degree $4$.
    \end{enumerate}    
    Then no two triangular faces incident with edges in $C$ have a common vertex in $\Int(C)$.
\end{lemma}
\textit{Proof.}
Suppose, for contradiction, that there are two triangular faces incident with edges in $C$ 
that share a common vertex in $\Int(C)$. 
Let these faces be $x_1x_2u$ and $y_1y_2u$, 
as depicted in Figure \ref{fig:CommonVertexTriangles}, 
where $u \in \Int(C)$ is the common vertex.

In Figure \ref{fig:CommonVertexTriangles}, $C_1$ and $C_2$ denote the parts of the cycle $C$:
$C_1$ is the path from $x_1$ to $y_1$ along $C$ (including $x_1$ and $y_1$); 
$C_2$ is the path from $x_2$ to $y_2$ along $C$ (including $x_2$ and $y_2$).
We consider three cases based on the degree of the vertex $u$.

\textsc{Case 1:} $\deg(u)=5$, and all vertices in $\Int(C)$ have degree $5$.

In this case, $u$ is connected to exactly five vertices. 
Since $u$ is connected to $x_1$, $x_2$, $y_1$, and $y_2$, 
it must have one additional neighbor. 
Let us denote it by $v$, 
and assume that it lies inside the cycle $C_2+u$, 
as shown in Figure \ref{fig:CommonVertexTriangles}(a).
This implies that $x_1y_1$ is an edge in $H$.
If $x_1y_1$ is an edge of the cycle $C$, then $\deg(x_1)=\deg(y_1)=3$.
This contradicts the assumption that all vertices in $C$ have degree at least $4$.
If $x_1y_1$ is a chord of $C$,
then this is contrary to Lemma \ref{lemma:The no chords}.

\textsc{Case 2:} $\deg(u)=5$, and 
at least one vertex in $\Int(C)$ has degree $6$.

In this case there is at most one vertex in $C$ of degree $5$.
If $\deg(x_1)=\deg(y_1)=4$, 
then the cycle $C_1+x_1y_1$, formed by $C_1$ together with the edge $x_1y_1$,
is a triangle $x_1y_1z$, where $z\in V(C_1)$.
This implies $\deg(z)=2$, which is impossible since all vertices in $C$ have degree at least $4$.

If $\deg(x_1)=4$ and $\deg(y_1)=5$, or vice versa, 
then all vertices of the cycle $C_1+x_1y_1$ have degree $4$, except for $x_1$ and $y_1$.
However, this contradicts Remark \ref{remark:case 3 of Lemma no chords}, 
which states that such a triangulation cannot exist.

\textsc{Case 3:} $\deg(u)=6$.

If both additional edges incident with $u$ lie inside the cycle $C_2+u$, 
then we can apply the same reasoning as in Cases 1 and 2, 
leading to a contradiction.

Let exactly one of the two additional edges incident with $u$ lie inside the cycle $C_2+u$
(see Figure \ref{fig:CommonVertexTriangles}(b)).
According to the conditions of the lemma, 
we can assume that $\deg(x_2)=\deg(y_2)=4$, 
and each vertex of $C_2$ also has degree $4$. 
Therefore, there exist vertices $v\in\Int(C)$ and $x_3,y_3\in V(C)$ such that 
$vx_2$, $vy_2$, $vx_3$, and $vy_3$ are edges in $H$.

If $\deg(v)=5$, we can repeat the reasoning from Cases 1 and 2, 
replacing $C_1$ with $C_2$ and $x_1$, $y_1$ with $x_3$, $y_3$.

If $\deg(v) = 6$, we proceed further 
(as shown in Figure \ref{fig:CommonVertexTriangles}(b)) to a vertex $w$.
Since by the condition of the lemma there are at most two vertices of degree $6$ inside $C$ 
and $\deg(u) = 6$ and $\deg(v) = 6$, it follows that $\deg(w) = 5$.
We can then apply the same reasoning as above, leading to a contradiction.

In all cases, we reach a contradiction.
Therefore, our initial assumption must be false, and the lemma is proved.

\phantom{1}

\textit{Proof} of Theorem \ref{Th:pure_triangulation}.
Suppose, on the contrary, that there exists a $3$-connected plane triangulation $G$ 
in which exactly one vertex $v$ has degree $n$ (where $n\neq5$)
and all other vertices have degree $5$.
Clearly, the $n$ neighbours of vertex $v$ form a cycle $C$ of length $n$.
Let $H=G-v$.
The graph $H$ is a plane triangulation of the $n$-cycle $C$,
and each vertex of $C$ has degree $4$ in the graph $H$.
Let $A\subset V(H)$ be the set of vertices of degree $5$.
Let $k=|A|$.
For convenience, we assume that $A\subset\Int(C)$.

\textsc{Step 1.}  
By Lemma \ref{lemma:on interior points}, $k=6$ (since $p_v=0$).
Thus, we have:
\begin{equation}\label{Cycle triangulation parameters}
    |V(H)|=n+6,\quad|E(H)|=2n+15,\quad|F(H)|=n+11.
\end{equation}

\textsc{Step 2.}
By virtue of Lemma \ref{lemma:The no chords}
no chord of the cycle $C$ is an edge of the graph $H$.

\textsc{Step 3.}
Due to step 2, we can conclude that each triangular face in $H$  has one, two, or three vertices from $A$.
Let $T_i$ denote the set of triangular faces with exactly $i$ vertices from $A$, where $i = 1, 2, 3$.
From Lemma \ref{lemma:no common vertex of triangles}, 
no two triangles in $T_1$ share a vertex from $A$. 
Therefore, the number of triangles in $T_1$ is at most $k = 6$.
On the other hand, by Lemma \ref{lemma:The no chords}, 
the edge of the triangle from $T_1$ opposite to the vertex from $A$ is an edge of the cycle $C$, 
and for each edge from $C$ there exists exactly one triangle from $T_1$ incident to this edge. 
This implies $|T_1|=n$.
Thus, $n\leq6$.
By Lemma \ref{lemma:interior_of_a_n-cycle (reinforced version)}, the cases $n=3$ and $n=4$ are impossible.
Additionally, we have assumed $n\neq5$. 
Thus, the only possibility is $n = 6$.

\begin{figure}[ht]
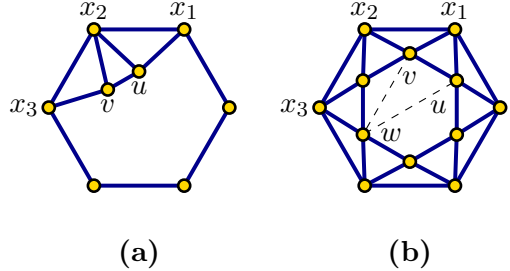

    \centering
    \ThreeTypesTriangle
    \caption[There are three types of triangular faces]
    {There are three types of triangular faces}
    \label{fig:ThreeTypesTriangle}
\end{figure}

\textsc{Step 4.}
Each vertex of $C$ is adjacent to exactly two vertices in $A$ 
(since its degree is $4$ and it is connected to two boundary vertices and two interior vertices). 
Therefore, each vertex of $C$ is part of exactly one triangle in $T_2$ 
(see Figure \ref{fig:ThreeTypesTriangle}(a)). 
Hence, $|T_2| = n = 6$.

Since
$$
|F(H)|=|T_1|+|T_2|+|T_3|+1=2n+1+|T_3|,
$$
it follows from (\ref{Cycle triangulation parameters}) that
$|T_3|=10-n=4$.

\textsc{Step 5.} 
The triangles of $T_1$ and $T_2$ are arranged as shown in Figure \ref{fig:ThreeTypesTriangle}(b).
All vertices of the triangles in $T_3$ are in $A$.
However, there cannot be an edge $uw$, otherwise $|\Int(x_1x_2x_3wu)|=2$;
there cannot be an edge $vw$, otherwise $|\Int(x_2x_3wu)|=1$.
In both cases we obtain a contradiction to Lemma
\ref{lemma:interior_of_a_n-cycle (reinforced version)}.
Therefore, $n\neq6$, which contradicts our earlier conclusion that $n=6$.

Hence, the theorem is completely proved.

\section{Triangulations with \texorpdfstring{$\Delta^+=1$}{dfc(G)=1}}
\label{section:Triangulations with dfc(G)=1}
Theorem \ref{Th:pure_triangulation} can be interpreted as follows.
If $G$ is a $3$-connected pentagulation of an $n$-gon and $n\neq5$, then $\Delta^+(G)>0$.
The purpose of this section is to prove the following statement.

\begin{theorem}\label{Th:dfc>1}
    Let $n$ be an integer and $n\geq3$, $n\neq5$.
    If $G$ is a $3$-connected pentagulation of an $n$-gon, then $\Delta^+(G)>1$.
\end{theorem}

Assume $\Delta^+(G)=1$ and let $p$ be the number of pentagonal faces of $G$. 
From (\ref{formula:DeficiencyDegreesGraph}), we obtain $p=n+8$, with
all vertices of $G$ have degree $3$ except for one, which has degree $4$.
Let $G^*$ be the dual of $G$.
In $G^*$, all vertices have degree $5$ except one, which has degree $n$,
and all faces are triangles except one, which is a square.
By Lemma \ref{lemma:Existence of chords}, the square has a chord.
The chord can either connect two vertices of degree $5$
or connect a vertex of degree $5$ to a vertex of degree $n$.
Now we add the chord of the square to $G^*$.
We get a new graph $H$.
This graph is a plane triangulation.
All but two or three of its vertices are of degree $5$.
If there are two exceptional vertices, they are adjacent.
One has degree $6$ and the other has degree $n+1$.
Among the three exceptional vertices there is one of degree $n$ and two adjacent vertices of degree $6$.

Conversely, suppose we have a plane triangulation $H$ containing one of the subgraphs
\begin{equation}\label{path one}
{\rm(i})\ \pathtwo{2pt}{4.5ex}{6}{6}\oplus\pathone{2pt}{n}\qquad\qquad
{\rm(ii})\ \pathtwo{2pt}{4.5ex}{n+1}{6}
\end{equation}
and all other vertices not in these subgraphs have degree $5$.
For the subgraph notation used here, see Remark \ref{remark:notation for subgraphs} below.
We get a graph $F$ containing exactly one quadrangular face after removing the single edge in these subgraphs. Therefore, the dual $F^*$ of $F$ is a pentagulation of an $n$-gon.
\begin{remark}
  If $n=4$, then condition \ref{path one}(ii) simply means that $H$ has exactly one vertex of degree $6$.
  If $n=6$, condition \ref{path one}(i) means that $H$ has three vertices of degree $6$,
  with two of them being adjacent.
\end{remark}

This reasoning shows that Theorem \ref{Th:dfc>1} is derived from the following statement
(see also Remark \ref{remark:About a theorem Grunbaum}).
\begin{theorem}\label{Th:No triangulations with dfc=1}
    Let $n$ be an integer and $n\geq3$, $n\neq5$.
    There is no $3$-connected plane triangulation containing one of the above subgraphs,
    and all vertices not contained in this subgraph are of degree $5$.
\end{theorem}
\textit{Proof.}
Let $G$ be a $3$-connected plane triangulation with subgraph (\ref{path one}i) or (\ref{path one}ii),
and all vertices of $G$ not contained in this subgraph are of degree $5$.
We first prove that if such a plane triangulation exists, then $n\leq8$.
We're going to have to examine a few cases.

\textsc{Case 1:} $G$ contains the  subgraph (\ref{path one}i).

If $v\in V(G)$ and $\deg(v)=n$, then $H=G-v$ is a plane triangulation of an $n$-gon $C$.
We may suppose  without loss of generality that $V(H)\subset\Int(C)\cup V(C)$.
Let $A=V(H)\cap\Int(C)$ and $|A|=k$. 
By Lemma \ref{lemma:on interior points}, we obtain $k=8$ (since $p_v=2$).

All vertices in $\Int(C)$ have degree $5$, and all vertices in $C$ have degree $4$,
except for the two adjacent vertices discussed below.
If both vertices of degree $6$ of $G$ are not neighbors of $v$,
then they are vertices of the same degree in $H$ and lie in $\Int(C)$.
If they are both neighbors of vertex $v$, then their degrees in $H$ are $5$
and they are adjacent and both lie in $C$.
Finally, if one of the vertices of degree $6$ is adjacent to $v$ and the other is not,
then there is exactly one vertex of degree $6$ in $\Int(C)$ and
there is exactly one vertex of degree $5$ in $C$.
Therefore,
Lemmas \ref{lemma:The no chords} and \ref{lemma:no common vertex of triangles} imply
$n\leq k=8$.

By Lemma \ref{lemma:interior_of_a_n-cycle (reinforced version)} the
case $n=3$ is impossible and by assumption $n\neq5$, so $n=4,6,7$, or $8$.

\begin{remark}
Note that the conditions of Lemma \ref{lemma:The no chords} are fulfilled in this case. 
Indeed, if there is at most one vertex of degree $5$ in $C$, 
then there is no requirement for internal vertices in condition (i) of Lemma \ref{lemma:The no chords}. 
If there are two vertices of degree $5$ in $C$, 
then all vertices inside $C$ also have degree $5$ 
which agrees with condition (ii) of Lemma \ref{lemma:The no chords}.
\end{remark}

\textsc{Case 2:} $G$ contains the  subgraph (\ref{path one}ii) and $n\neq4$.
Let $v\in V(G)$ and $\deg(v)=n+1$.
Then $H=G-v$ is a plane triangulation of an $(n+1)$-gon $C$.
Let $A$ be the set of all internal vertices of the cycle $C$.
Denote $k=|A|$. By Lemma \ref{lemma:on interior points}, we obtain $k=7$ (since $p_v=1$).

Lemmas \ref{lemma:The no chords}(i) and \ref{lemma:no common vertex of triangles}(i) imply
$n+1\leq7$, so $n\leq6$.
By Lemma \ref{lemma:interior_of_a_n-cycle (reinforced version)}.2 the
case $n+1=4$ ($n=3$) is impossible and by assumption $n\neq5$ and $n\neq4$, 
so $n=6$.

If $n=6$, then we have a triangulation of the $7$-gon $C$. 
Each vertex of $C$ has degree $4$, except for one vertex of degree $5$. 
In $\Int(C)$ there are $7$ vertices, each of which has degree $5$.
\begin{figure}[ht]
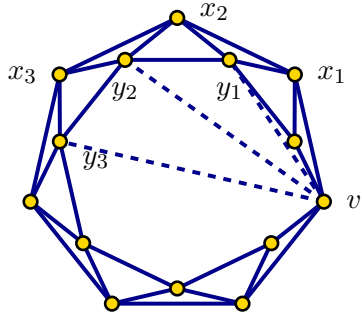

    \centering
    \HeptagonWithSevenNodeInside{2}{0.09}
    \caption[A heptagon with seven vertices inside]
    {A heptagon with seven vertices inside}
    \label{fig:HeptagonWithSevenNodeInside}
\end{figure}
In Figure \ref{fig:HeptagonWithSevenNodeInside}, vertex $v$ has degree $5$. 
The other vertices of the cycle $C$ have degree $4$. 
The second vertex of the additional edge incident to vertex $v$ can be one of three vertices: 
$y_1$, $y_2$, or $y_3$, by symmetry.

However, none of these cases is possible.
By Corollary 
\ref{corollary:interior_of_a_n-cycle (weak variant)}:
if $vy_1$ is an edge, then there is only one vertex inside the triangle $vx_1y_1$;
if $vy_2$ is an edge, then there is only two vertices inside the square  $vx_1x_2y_2$;
if $vy_3$ is an edge, then there is only three vertices inside the pentagon $vx_1x_2x_3y_3$.

\textsc{Case 3:} $G$ contains the  subgraph (\ref{path one}ii) and $n=4$.

Then $G$ contains 
exactly one vertex of degree $6$, say $v$, while all other vertices have degree $5$.
All neighbours of $v$ form a cycle $C$ of length $6$. 
Let $A$ be the set of internal vertices of $C$ and $k=|A|$.
By Lemma \ref{lemma:on interior points}, we obtain that $k=6$ (since $p_v$=0).
As we seen in Step 5 of Theorem \ref{Th:pure_triangulation}, 
the existence of a triangulation of a $6$-gon 
whose vertices all have degree $4$ 
and whose interior vertices have degree $5$ leads to a contradiction.
Therefore, this case is impossible.

\textsc{Case 4:} $G$ contains the  subgraph (\ref{path one}i) and $n=4$.

We investigate this case, as well as the cases when $n=6$, $7$, and $8$, 
using the program \texttt{plantri} \cite{Plantri}.
We use \texttt{plantri} together with the plug-in \texttt{allowed\_deg},
which is distributed with \texttt{plantri}, and
allows us to specify which degrees are permitted.

As we have seen above, the number of vertices of the desired planar triangulation $G$
must be equal to $n+8+1=n+9$. Because $n=4$, we get $|V(G)|=13$.   
We must compute planar triangulations of the following type:
Those that contain two vertices of degree $6$ and one vertex of degree $4$ in addition to vertices of degree $5$.
We can use the following command-line parameters:
\begin{align*}
    \texttt{plantri\_ad 13 -a -F5F6\_2\textasciicircum2F4\_1\textasciicircum1}
\end{align*}
which generates the $3$-connected triangulations with $13$ vertices having
the required number of vertices of degree $6$ and degree $4$.
The switch \texttt{-a}
indicates that the graphs will be output in ASCII format.
In this case, only one triangulation was found in which two vertices of degree $6$ are not adjacent.

\textsc{Case 5:} $G$ contains the  subgraph (\ref{path one}i) and $n=6$.

The number of vertices of the desired planar triangulation
must be equal to $n+9=15$:
\begin{align*}
  \texttt{plantri\_ad 15 -a -F5F6\_3\textasciicircum3}
\end{align*}
Only one triangulation was found that has three vertices of degree $6$.
However, there is no pair of adjacent vertices of degree $6$ in the found triangulation.

\textsc{Case 6:} $G$ contains the  subgraph (\ref{path one}i) and $n=7$ or $n=8$.
The number of vertices of the desired planar triangulation must be equal 
to $16$, if $n=7$, and to $17$, if $n=8$:
\begin{align*}
  \texttt{plantri\_ad 16 -a -F5F6\_2\textasciicircum2F7\_1\textasciicircum1} \\  
  \texttt{plantri\_ad 17 -a -F5F6\_2\textasciicircum2F8\_1\textasciicircum1} 
\end{align*}
In both cases, the output is zero.

The theorem is completely proved.
\begin{remark}
    A command string of the form 
    \begin{align*}
      \texttt{plantri\_ad 16 -a -F5F6\_2\textasciicircum2F7\_1\textasciicircum1}       
    \end{align*}    
    outputs only those triangulations whose vertex degrees are $5$, $6$, or $7$.
    In this case, only those with exactly two vertices of degree $6$ and exactly one vertex of degree $7$
    are generated.
    For more details, please refer to the manual \cite{Plantri}.
\end{remark}

\section{Triangulations with \texorpdfstring{$\Delta^+=2$}{dfc(G)=2}}
\label{section:Triangulations with dfc(G)=2}
We assume in this section that $n\neq5$.
If $\Delta^+(G)=2$, it follows from (\ref{formula:DeficiencyDegreesGraph}) that $p=n+10$.
For the graph $G$, there are two possibilities:
\begin{enumerate}[label=$\Delta$.2.\Roman*,align=left,left=\the\parindent,labelsep=-0.4em]
      \item\label{item:dfc2 degree 5}
      All vertices have degree $3$ except one, which has degree $5$;
      \item\label{item:dfc2 degree 4 4}
      Two vertices are degree $4$, all others degree $3$.
\end{enumerate}
Let $G^*$ be the dual of $G$.
Exactly one vertex of $G^*$ has degree $n$, and all other vertices have degree $5$.

In case \ref{item:dfc2 degree 5}, the graph $G^*$ has exactly one pentagonal face, 
and all the other faces are triangles.
By Lemma \ref{lemma:Existence of chords}, the pentagon has two chords.
Let one chord join the vertices $u$ and $v$, and let the other chord join the vertices $u$ and $w$.
There are three possibilities:
\begin{enumerate}[label=(\roman*),align=left,left=\the\parindent]
      \item $\deg(u)=\deg(v)=\deg(w)=5$;
      \item $\deg(u)=n$, $\deg(v)=\deg(w)=5$;
      \item $\deg(u)=\deg(v)=5$ and $\deg(w)=n$ or $\deg(u)=\deg(w)=5,\deg(v)=n$.
\end{enumerate}
We add the two chords of the pentagon to $G^*$.
As a result, we obtain a plane triangulation $H$
that contains one of the following subgraphs:
 $$
 {\rm(i)}\ \paththree{2pt}{5ex}{6}{7}{6}\oplus\pathone{2pt}{n}\phantom{mmmm}
 {\rm(ii)}\ \paththree{2pt}{5ex}{6}{n+2}{6}\phantom{mmmm}
 {\rm(iii)}\ \paththree{2pt}{5ex}{6}{7}{n+1}
 $$
All other vertices of $H$ have degree $5$.
\begin{remark}\label{remark:notation for subgraphs}
Here we will use the following notation for subgraphs of a graph $G$.
This is how we denote three subgraphs: exactly one vertex, a path of length $1$ (edge), and a path of length $2$:
$$
\pathone{2pt}{n} \phantom{mmmm}
\pathtwo{2pt}{5ex}{6}{6} \phantom{mmmm}
\paththree{2pt}{5ex}{6}{7}{6}
$$
The number denotes the degree of the corresponding vertex in the graph $G$.
In the following, we will also encounter the following subgraphs:
$$
  \pathfour{2pt}{5ex}{6}{7}{7}{6}\qquad
  \raisebox{-2.5ex}{
  \xstar{2pt}{5ex}{6}{8}{6}{6}
  }\qquad
  \raisebox{-3.3ex}{
  \xtriangle{2pt}{5ex}{7}{7}{7}
  }
$$
If $X$ and $Y$ are two graphs, then
$X\oplus Y$ denotes their disjoint union.
\end{remark}
\begin{remark}
    If $n=3$, then condition (ii) for the graph $H$ means that the two vertices of degree $6$ have a common neighbor of degree $5$.
    If $n=4$, then condition (iii) in $H$ simplifies to
    $\pathtwo{2pt}{5ex}{6}{8}$.
\end{remark}
In case \ref{item:dfc2 degree 4 4}, the graph $G^*$ has exactly two square faces in addition to the triangular faces.
According to Lemma \ref{lemma:Existence of chords}, each of these squares contains a chord.
Let $u$, $v$ be the vertices of one chord, and $x$, $y$ be the vertices of the other chord.

If we add the chords to $G^*$, we get a plane triangulation $H$ in which,
in addition to the vertices of degree $5$, there are three to five vertices of degrees
$n$, $n+1$, $n+2$, $6$, or $7$.
depending on the location of the two square faces in $G^*$.
If $\{u,v\}\cap\{x,y\}=\varnothing$, then $H$ has one of the following subgraphs:
$$
{\rm(i)}\ \pathtwo{2pt}{5ex}{6}{6}\oplus\pathtwo{2pt}{5ex}{6}{6}\oplus\pathone{2pt}{n}\phantom{mmmm}
{\rm(iii)}\ \pathtwo{2pt}{5ex}{6}{6}\oplus\pathtwo{2pt}{5ex}{6}{n+1}
$$
If $\{u,v\}\cap\{x,y\}\neq\varnothing$ (e.g. $v=x$), then $H$ has one of the following subgraphs:
$$
{\rm(ii)}\ \paththree{2pt}{5ex}{6}{7}{6}\oplus\pathone{2pt}{n}\phantom{mmmm}
{\rm(iv)}\ \paththree{2pt}{5ex}{6}{7}{n+1}\phantom{mmmm}
{\rm(v)}\ \paththree{2pt}{5ex}{6}{n+2}{6}
$$
\begin{remark}
    As mentioned earlier, some of these five conditions simplify for $n=3$ or $n=4$.
\end{remark}
By combining cases \ref{item:dfc2 degree 5} and \ref{item:dfc2 degree 4 4}, we can now state the following.
\begin{theorem}\label{Th:pentagulation_eqiv_triangulation_n+10}
    Let $n\geq3$ and $n\neq5$.
    A pentagulation of an $n$-gon in which the number of pentagons is $p=n+10$ exists,
    if and only if
    there exists a plane triangulation in which there exists one of the following subgraphs:
    \begin{align*}
        {\rm(i)} & \pathtwo{2pt}{5ex}{6}{6}\oplus\pathtwo{2pt}{5ex}{6}{6}\oplus\pathone{2pt}{n}\phantom{mmmm} \\
        {\rm(ii)} & \paththree{2pt}{5ex}{6}{7}{6}\oplus\pathone{2pt}{n}\phantom{mmmm}\\
        {\rm(iii)}& \pathtwo{2pt}{5ex}{6}{6}\oplus\pathtwo{2pt}{5ex}{6}{n+1}\\
        {\rm(iv)}& \paththree{2pt}{5ex}{6}{7}{n+1}\\
        {\rm(v)}& \paththree{2pt}{5ex}{6}{n+2}{6}\phantom{mmmm}
    \end{align*}
    and all other vertices have degree $5$.
\end{theorem}

\textit{Proof.}
The necessity was proved before the formulation of the theorem.

Conversely, suppose we have a plane triangulation $H$ in which there is one of
the indicated subgraphs and all other vertices of $H$ have degree $5$.
A graph $F$ is obtained by removing two edges from $H$, 
namely those forming subgraphs (i)--(v), if present in $H$.
As a result, the graph $F$ has all but one vertex of degree $5$.
The exceptional vertex of $F$ has degree $n$.
Furthermore, all but one or two faces of the graph $F$ are triangles.
If the edges to be removed were part of the boundary of a triangular face
of the graph $H$, then the graph $F$ contains exactly one pentagonal face.
Otherwise, the graph $F$ will contain exactly two square faces.

It follows that all but one of the faces of the dual $F^*$ are pentagons.
The exceptional face is an $n$-gon.
Furthermore, all vertices of $F^*$ are of degree $3$, 
except for either one vertex of degree $5$ or two vertices of degree $4$.
Consequently, $F^*$ is a pentagulation of an $n$-gon with $p=n+10$ pentagons.
This completes the proof of the theorem.

\begin{theorem}\label{Th:Triangulations with dfc=2}
    If $n=4,6$, or $8$, then there exists a $3$-connected plane triangulation containing one of the above subgraphs, and all vertices not contained in this subgraph are of degree $5$.
    If $n=3,7$ or $9\leq n\leq13$, then such triangulations do not exist.

    Furthermore,
    for each $n=4,5,6$ there is a unique (up to isomorphism) $3$-connected pentagulation with $\Delta^+=2$.
\end{theorem}
\textit{Proof.}
As we know the number of vertices in the sought triangulations is $n+11$.
To find the specified triangulations, we will use \texttt{plantri}.
Table \ref{tabl: degree sequences of triangulations (dfc=2)} 
lists the degree sequences of triangulations corresponding to the types specified in Theorem \ref{Th:pentagulation_eqiv_triangulation_n+10}.
\begin{table}[ht!]
\begin{center}
\small
\begin{tabular}{|l|l|l|l|l|}
  \hline
       &          3  &     4       &      6      &  $n\geq7$      \\
  \hline
  (i)  & (3,6,6,6,6) & (4,6,6,6,6) & (6,6,6,6,6) & $(6,6,6,6,n)$  \\
  (ii) & (3,6,6,7)   & (4,6,6,7)   & (6,6,6,7)   & $(6,6,7,n)$    \\
  (iii)& (4,6,6,6)   & (6,6,6)     & as in (ii)  & $(6,6,6,n+1)$  \\
  (iv) & (4,6,7)     & (6,7)       & (6,7,7)     & $(6,7,n+1)$    \\
  (v)  & (6,6,8)     & as in (iii) & (6,6,8)     & $(6,6,n+2)$    \\
  \hline
\end{tabular}
\caption{Degree sequences of triangulations. Case $\Delta^+ = 2$.}
\label{tabl: degree sequences of triangulations (dfc=2)}
\end{center}
\end{table}
\begin{remark}
Only degrees other than $5$ are shown here and below. 
Therefore, the necessary number of fives must be added to each sequence 
to represent the degrees of all vertices in the graph, that is, to obtain a 
degree sequence of the graph.
\end{remark}
\begin{table}[ht!]
\begin{center}
\small
\begin{tabular}{|l|l|l|l|l|l|l|}
  \hline
       & 3 & 4 & 6 & 7 & 8 & $n\geq9$ \\
  \hline
  (i)  & 0 & 3 & 3 & 2 & 1 & 0\\
  (ii) & 1 & 2 & 1 & 3 & 1 & 0\\
  (iii)& 1 & 1 & --& 0 & 0 & 0\\
  (iv) & 0 & 0 & 0 & 0 & 0 & 0\\
  (v)  & 1 & --& 0 & 0 & 0 & 0\\
  \hline
\end{tabular}
\caption{Number of triangulations returned by \texttt{plantri}. Case $\Delta^+ = 2$.}
\label{tabl: number of triangulations (dfc=2)}
\end{center}
\end{table}
With \texttt{plantri} we can obtain all triangulations with $n+11$ vertices,
having the degree sequences given in Table \ref{tabl: degree sequences of triangulations (dfc=2)}.
The number of triangulations of each type, as obtained with \texttt{plantri} is shown in Table
\ref{tabl: number of triangulations (dfc=2)}.
For example,
\begin{align*}
\texttt{plantri\_ad 18 -a -F5F6\_2\textasciicircum2F7\_2\textasciicircum2}
\end{align*}
returns $3$, as shown in the second row and fourth column of Table
\ref{tabl: number of triangulations (dfc=2)}.
A detailed analysis shows that for $n=3$ and $n=7$, none of the graphs contains the subgraphs listed in Theorem
\ref{Th:Triangulations with dfc=2}.

\begin{figure}[ht]
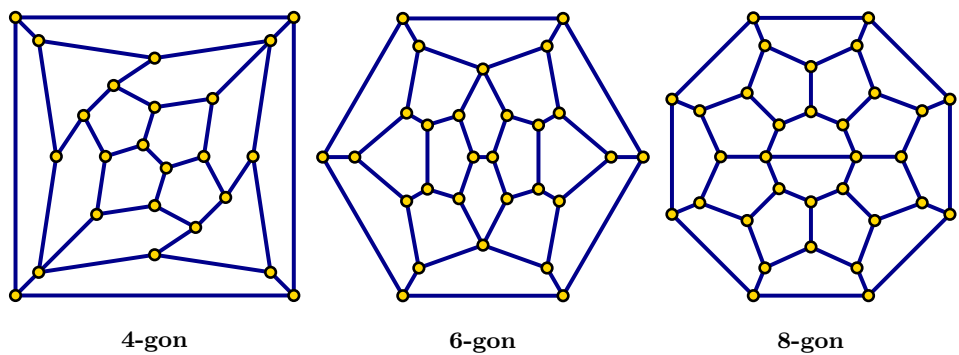

    \centering
    \SquareMP{21.8mm}{0.7mm}{\footnotesize 4-gon}
    \HexagonMP{21.3mm}{0.7mm}{\footnotesize 6-gon}
    \OktagonMP{20mm}{0.7mm}{\footnotesize 8-gon}
    \caption[Minimal  $3$-connected pentagulations of 4,6,8-gons]
    {Minimal  $3$-connected pentagulations of 4,6,8-gons.}
    \label{fig:SquareHexagonOktagonMP}
\end{figure}

\begin{figure}[ht]
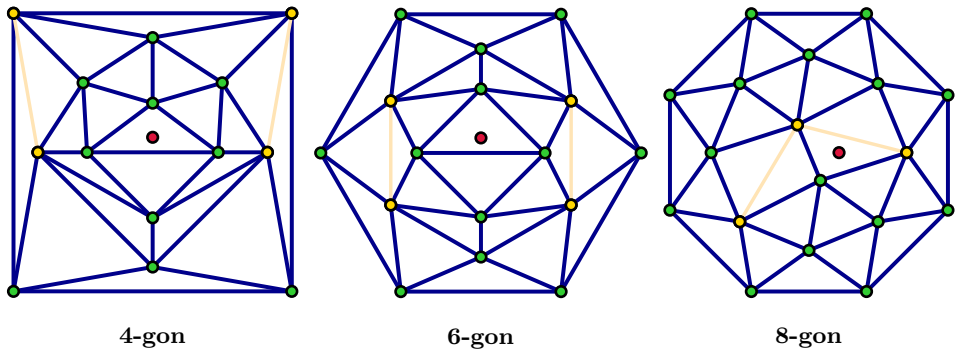

    \centering
    \SquareT{21.8mm}{0.7mm}{\footnotesize 4-gon}
    \HexagonT{21.3mm}{0.7mm}{\footnotesize 6-gon}
    \OktagonT{20mm}{0.7mm}{\footnotesize 8-gon}
    \caption[Triangulations that correspond to\\ minimal $3$-connected pentagulations of 4,6,8-gons]
    {Triangulations that correspond to\\ minimal $3$-connected pentagulations of 4,6,8-gons.}
    \label{fig:SquareHexagonOktagonT}
\end{figure}

For $n=4,6,8$, among the triangulations found, there are triangulations that contain subgraphs from Theorem
\ref{Th:Triangulations with dfc=2}.
Each of these triangulations corresponds to exactly one pentagulation for  $n = 4,6,8$.
Figure \ref{fig:SquareHexagonOktagonMP} depicts these pentagulations and
Figure \ref{fig:SquareHexagonOktagonT} illustrates the corresponding plane triangulations.
The theorem is proved.

\begin{remark}\label{remark:colors and lieders}
    We make a few observations about Figure \ref{fig:SquareHexagonOktagonT}.

    First, not all edges and vertices of each triangulation are displayed in the figure.    
    Exactly one vertex, referred to as the leading vertex or simply the leader, 
    is omitted in the figure.
    To form a plane triangulation, a leader vertex must be introduced outside the 
    $n$-gon and connected to every vertex of the $n$-gon with edges.
    The red vertex is a visual reminder of the leader's existence and is not part of the configuration.    
    We believe this approach enhances the symmetry and visual appeal of the drawings.

    Second, some edges and the vertices incident with them in each triangulation are highlighted in golden color.
    Removing all such golden edges results in a graph in which all vertices have degree $5$, 
    except for one vertex with degree $n$.
    When calculating the degree of any vertex, it is important to include the leader and all edges incident with it.
    The dual of this graph corresponds to the desired pentagulation.
    
    Third, selecting a different vertex as the leader results in 
    an alternative visual representation of the same graph.    
    For instance, Figure \ref{fig:SquareTMPF} illustrates a triangulation in which
    a vertex of degree $5$ is chosen as the leader.
    The dual of the graph obtained by removing two golden edges from this triangulation
    represents a pentagulation of a square.
    The two pentagulations of the square shown in Figures \ref{fig:SquareHexagonOktagonT}($4$-gon) and \ref{fig:SquareTMPF}(b) are isomorphic.
\begin{figure}[ht]
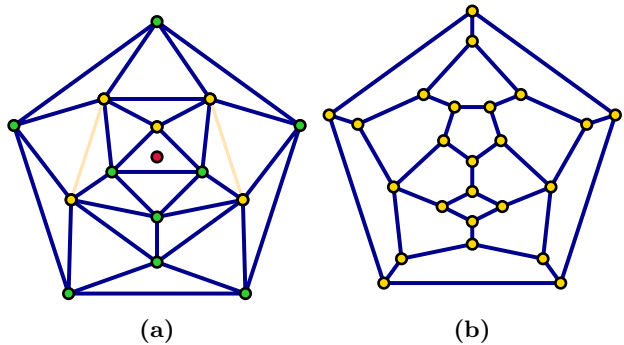

    \centering
    \SquareTF{20mm}{0.7mm}{\footnotesize (a)}
    \SquareMPF{20mm}{0.7mm}{\footnotesize (b)}
    \caption[A triangulation of a square with a leading vertex of degree $5$]
    {A triangulation of a square with a leading vertex of degree $5$ and \\its corresponding pentagulation.}
    \label{fig:SquareTMPF}
\end{figure}
\end{remark}

\section{Triangulations with \texorpdfstring{$\Delta^+=3$}{dfc(G)=3}}
\label{section:Triangulations with dfc(G)=3}
We maintain the assumption that $n\neq5$.
If $G$ is a $3$-connected pentagulation of an $n$-gon with $\Delta^+(G) = 3$,
then due to (\ref{formula:DeficiencyDegreesGraph}) we obtain $p=n+12$.
We distinguish three cases:
\begin{enumerate}[label=$\Delta$.3.\Roman*,align=left,left=\the\parindent,labelsep=-0.4em]
      \item\label{item:one vertex of degree 6}
      All vertices of $G$ have degree $3$ except for one, which has degree $6$;
      \item\label{item:two vertices of degrees 4 and 5}
      Two vertices of $G$ are of degrees $4$ and $5$. All others are of degree $3$;
      \item \label{item:three vertices of degree 4}
      Three vertices of $G$ have degree $4$, while all others have degree $3$.
\end{enumerate}
The dual $G^*$ of $G$ has the following properties.
All vertices of $G^*$ have degree $5$ except for one with degree $n$.
All faces of $G^*$ are triangles except at most three, which are hexagons, pentagons, or squares.
In case \ref{item:one vertex of degree 6}, $G^*$  has one hexagonal face;
in case \ref{item:two vertices of degrees 4 and 5}, one pentagonal and one quadrangular face; 
and 
in case \ref{item:three vertices of degree 4}, exactly three quadrangular faces.

In view of Lemma \ref{lemma:Existence of chords},
there are three chords in the non-triangular faces that, collectively, 
subdivide these faces into four triangular faces.
The graphs with three edges and at most one isolated vertex are listed as follows:
\begin{align*}
 1.\quad & P_3 & & P_3\oplus K_1\\
 2.\quad & P_2\oplus K_2 & & P_2\oplus K_2\oplus K_1\\
 3.\quad & K_2\oplus K_2\oplus K_2 & & K_2\oplus K_2\oplus K_2\oplus K_1\\
 4.\quad & K_{1,3} & & K_{1,3}\oplus K_1\\
 5.\quad & K_3 & & K_3\oplus K_1
\end{align*}
Standard notation is used:
$\oplus$ indicates the disjoint union of graphs;
$P_s$ denotes a path with $s$ edges;
$K_s$ is a complete graph on $s$ vertices;
and $K_{s,t}$ represents a complete bipartite graph with partitions of size $s$ and $t$.
The graph $K_{1,3}$ is commonly referred to as the claw.

We now add the three chords to $G^*$.
This produces a plane triangulation $H$ that includes,
in addition to the vertices of degree $5$,
three to seven vertices with degrees $n,n+1,n+2,n+3$, $6$, $7$, or $8$.
The graph $H$ contains a subgraph of one of the following types:
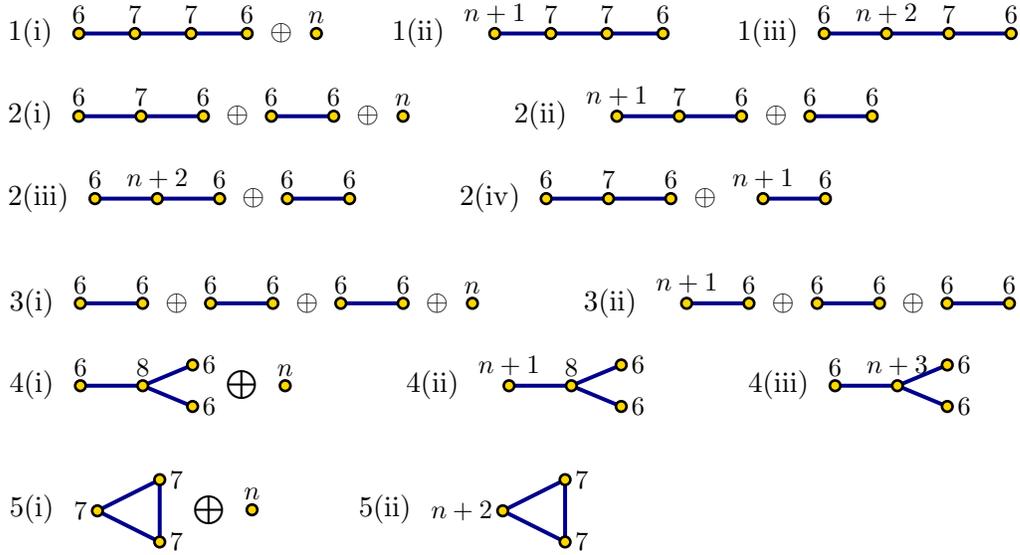
\begin{align*}
&
1({\rm i})\ \pathfour{2pt}{4.5ex}{6}{7}{7}{6}\oplus\pathone{2pt}{n}\qquad
1({\rm ii})\ \pathfour{2pt}{4.5ex}{n+1}{7}{7}{6}\qquad
1({\rm iii})\ \pathfour{2pt}{5ex}{6}{n+2}{7}{6}
\vphantom{\int_I^I}
\\
&
2({\rm i})\ \paththree{2pt}{5ex}{6}{7}{6}\oplus\pathtwo{2pt}{5ex}{6}{6}\oplus\pathone{2pt}{n}\phantom{MMM}
2({\rm ii})\ \paththree{2pt}{5ex}{n+1}{7}{6}\oplus\pathtwo{2pt}{5ex}{6}{6}
\vphantom{\int_I^I}
\\
&
2({\rm iii})\ \paththree{2pt}{5ex}{6}{n+2}{6}\oplus\pathtwo{2pt}{5ex}{6}{6}\phantom{MMM}
2({\rm iv})\ \paththree{2pt}{5ex}{6}{7}{6}\oplus\pathtwo{2pt}{5ex}{n+1}{6}
\end{align*}
\begin{align*}
&
3({\rm i})\ \pathtwo{2pt}{5ex}{6}{6}\oplus\pathtwo{2pt}{5ex}{6}{6}\oplus\pathtwo{2pt}{5ex}{6}{6}\oplus\pathone{2pt}{n}\phantom{MMM}
3({\rm ii})\ \pathtwo{2pt}{5ex}{n+1}{6}\oplus\pathtwo{2pt}{5ex}{6}{6}\oplus\pathtwo{2pt}{5ex}{6}{6}
\vphantom{\int_I^I}
\\
&
\raisebox{2.6ex}{$4$(i)}\ \xstar{2pt}{5ex}{6}{8}{6}{6}\raisebox{2.6ex}{$\bigoplus$\ \pathone{2pt}{n}}\phantom{MMM}
\raisebox{2.6ex}{$4$(ii)}\ \xstar{2pt}{5ex}{n+1}{8}{6}{6}\phantom{MMM}
\raisebox{2.6ex}{$4$(iii)}\ \xstar{2pt}{5ex}{6}{n+3}{6}{6}
\vphantom{\int_I^I}
\\
&
\raisebox{3.5ex}{$5$(i)}\ \xtriangle{2pt}{5ex}{7}{7}{7}\raisebox{3.5ex}{$\bigoplus$\ \pathone{2pt}{n}\phantom{MMM}}
\raisebox{3.5ex}{$5$(ii)}\ \xtriangle{2pt}{5ex}{7}{7}{n+2}
\end{align*}

\begin{theorem}\label{Th:pentagulation_eqiv_triangulation_n+12}
    Let $n\geq3$ and $n\neq5$.
    A pentagulation of an $n$-gon with $p=n+12$ pentagons exists
    if and only if
    there is a plane triangulation that contains a subgraph of one of the above $14$ types.
    All other vertices not in this subgraph have degree $5$.
\end{theorem}
The proof of this theorem is similar to the proof of Theorem \ref{Th:pentagulation_eqiv_triangulation_n+10}
and we omit it.

We exclude the cases $n=4,6,8$ from the following theorem, 
as our aim is to identify minimal pentagulations of the $n$-gon.
\begin{theorem}\label{Th:Triangulations with dfc=3}
    If $n=3,7,9,10,11,12$, then there is a $3$-connected plane triangulation 
    containing one of the above subgraphs, and all vertices not contained in this subgraph are of degree $5$.
    No such triangulations exist for $n=13$.

    Furthermore,
    if $\Delta^+=3$, then for $n=3,7,10$ there are exactly three
    $3$-connected pentagulations each, up to isomorphism;
    for $n=9$ there are four such pentagulations;
    and for $n=11$ and $12$ there is exactly one unique pentagulation.
\end{theorem}
\textit{Proof.}
The number of vertices in the desired triangulations is $n+13$.
For each $n$, the subgraphs 2(ii) and 2(iv) have the same degree sequences;
therefore, the row corresponding to subgraph 2(iv) is omitted
in Table \ref{tabl: degree sequences of triangulations (dfc=3)}.
For $n=3$, subgraphs 2(iii) and 4(iii)
have the same degree sequences, as do subgraphs 2(ii) and 4(i) for $n=7$.

\begin{table}[ht!]
\begin{center}
\small
\begin{tabular}{|l|l|l|l|}
  \hline
  &           3          &      7        &   $n\geq9$      \\
  \hline
  1(i)  & (3,6,6,7,7)   &(6,6,7,7,7)    & $(6,6,7,7,n)$    \\
  1(ii) & (4,6,7,7)     &(6,7,7,8)      & $(6,7,7,n+1)$    \\
  1(iii)& (6,6,7)       &(6,6,7,9)      & $(6,6,7,n+2)$    \\
  2(i)  &(3,6,6,6,6,7)  &(6,6,6,6,7,7)  & $(6,6,6,6,7,n)$  \\
  2(ii) & (4,6,6,6,7)   &(6,6,6,7,8)    & $(6,6,6,7,n+1)$  \\
  2(iii)& (6,6,6,6)     &(6,6,6,6,9)    & $(6,6,6,6,n+2)$  \\
  3(i)  &(3,6,6,6,6,6,6)&(6,6,6,6,6,6,7)& $(6,6,6,6,6,6,n)$\\
  3(ii) &(4,6,6,6,6,6)  &(6,6,6,6,6,8)  & $(6,6,6,6,6,n+1)$\\
  4(i)  &(3,6,6,6,8)    &as in 2(ii)    & $(6,6,6,8,n)$    \\
  4(ii) &(4,6,6,8)      &(6,6,8,8)      & $(6,6,8,n+1)$    \\
  4(iii)&as in 2(iii)   &(6,6,6,10)     & $(6,6,6,n+3)$    \\
  5(i)  &(3,7,7,7)      &(7,7,7,7)      & $(7,7,7,n)$      \\
  5(ii) &(7,7)          &(7,7,9)        & $(7,7,n+2)$      \\
  \hline
\end{tabular}
\caption{Degree sequences of triangulations. Case $\Delta^+=3$.}
\label{tabl: degree sequences of triangulations (dfc=3)}
\end{center}
\end{table}

We use \texttt{plantri} to find the specified triangulations.
Table \ref{tabl: degree sequences of triangulations (dfc=3)}
shows all possible degree sequences corresponding to the subgraphs listed above.
Table \ref{tabl: number of triangulations (dfc=3)} provides the number of triangulations 
of each type obtained using \texttt{plantri}.
The numbers in brackets in this table indicate the number of triangulations that contain a subgraph of the corresponding type.
A blank space indicates the absence of triangulations with the corresponding degree sequence.
Some numbers in the table appear in bold coloured text.
Numbers corresponding to the same degree sequences are highlighted in the same colour.
\newcommand\bfred[1]{\textcolor{red}{\bf#1}}
\newcommand\bfblue[1]{\textcolor{blue}{\bf#1}}
\begin{table}[ht!]
\begin{center}
\small
\begin{tabular}{|l|l|l|l|l|l|l|l|}
  \hline
       &    3         & 7               & 9    & 10   & 11   & 12   & 13 \\
  \hline
1(i)   & 5(1)         & 5(1)            & 7(3) & 6(6) & 3(3) & 6(2) & 5(1) \\
1(ii)  & 2(0)         &                 &      &      &      &      &      \\
1(iii) &              &                 &      &      &      &      &      \\
2(i)   & 8(4)         & 23(5)           & 8(7) &      & 2(2) & 2(2) & 3(1) \\
2(ii)  & \bfblue{4}(2)& \bfred{6}(0)    &      &      &      &      &      \\
2(iii) & \bfred{2}(1) &                 &      &      &      &      & 1(1) \\
2(iv)  & \bfblue{4}(0)& \bfred{6}(3)    &      &      &      &      &      \\
3(i)   & 3(2)         & 16(7)           & 2(2) &      &      & 1(1) &      \\
3(ii)  & 5(3)         & 2(0)            &      &      &      &      &      \\
4(i)   & 4(0)         & \bfred{6}(0)    & 7(3) & 6(6) & 1(1) &      &      \\
4(ii)  &              &                 &      &      &      &      &      \\
4(iii) & \bfred{2}(0) & 2(0)            &      &      &      &      & 1(0) \\
5(i)   & 2(1)         &                 &      & 3(3) &      &      &      \\
5(ii)  & 1(0)         &                 &      &      &      &      &      \\
  \hline
\end{tabular}
\caption{Number of triangulations returned by \texttt{plantri}. Case $\Delta^+=3$.}
\label{tabl: number of triangulations (dfc=3)}
\end{center}
\end{table}
For each triangulation found, all subgraphs of the corresponding type were computed.
A single triangulation may contain more than one subgraph of the corresponding type.
For instance, only two triangulations exist for the degree sequence $(6,6,6,6)$:
\begin{align*}
    \texttt{1) bcdef,afghc,abhijd,acjke,adklmf,aemgb,bfmnoh,bgoic,chopj,cipkd,}\\
    \texttt{djple,ekpnm,elngf,gmlpo,gnpih,ionlkj;} \\
    \texttt{2) bcdef,afghc,abhid,acijke,adklmf,aemgb,bfmnh,bgnoic,chojd,diopk,}\\
    \texttt{djple,ekpnm,elngf,gmlpoh,hnpji,jonlk.}
\end{align*}
(Recall that fives are omitted in degree sequences of triangulations, 
as their count is determined by Euler's theorem.)
Each graph here consists of $16$ vertices.
The vertices are labeled using letters \texttt{a, b, \ldots, p}.
For each vertex, its adjacent vertices are listed.
For example, vertex \texttt{a} is adjacent to \texttt{b, c, d, e, f} in both graphs.
The first graph contains no subgraphs of types 2(iii) or 4(iii).
The second graph lacks subgraphs of type 4(iii), but includes several of type 2(iii).
This graph, in fact, has $4$ vertices of degree $6$: \texttt{d, e, h, n}, 
while all others have degree $5$. 
There are four subgraphs of type 2(iii) in this graph:
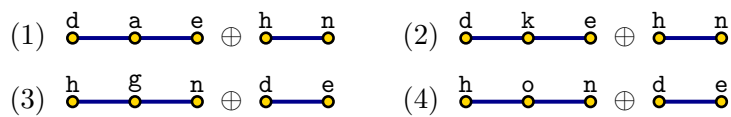
\begin{align*}
&
(1)\
\paththree{2pt}{5ex}{\texttt{d}}{\texttt{a}}{\texttt{e}}\oplus\pathtwo{2pt}{5ex}{\texttt{h}}{\texttt{n}}
\qquad
(2)\
\paththree{2pt}{5ex}{\texttt{d}}{\texttt{k}}{\texttt{e}}\oplus\pathtwo{2pt}{5ex}{\texttt{h}}{\texttt{n}}\\
&
(3)\
\paththree{2pt}{5ex}{\texttt{h}}{\texttt{g}}{\texttt{n}}\oplus\pathtwo{2pt}{5ex}{\texttt{d}}{\texttt{e}}
\qquad
(4)\
 \paththree{2pt}{5ex}{\texttt{h}}{\texttt{o}}{\texttt{n}}\oplus\pathtwo{2pt}{5ex}{\texttt{d}}{\texttt{e}}
\end{align*}

From each triangulation containing a subgraph of one of the 14 types, we remove three edges of the subgraph.
We then construct the dual of the resulting graph, which forms a pentagulation of the corresponding $n$-gon.
This process results in the following.

\begin{figure}[ht]
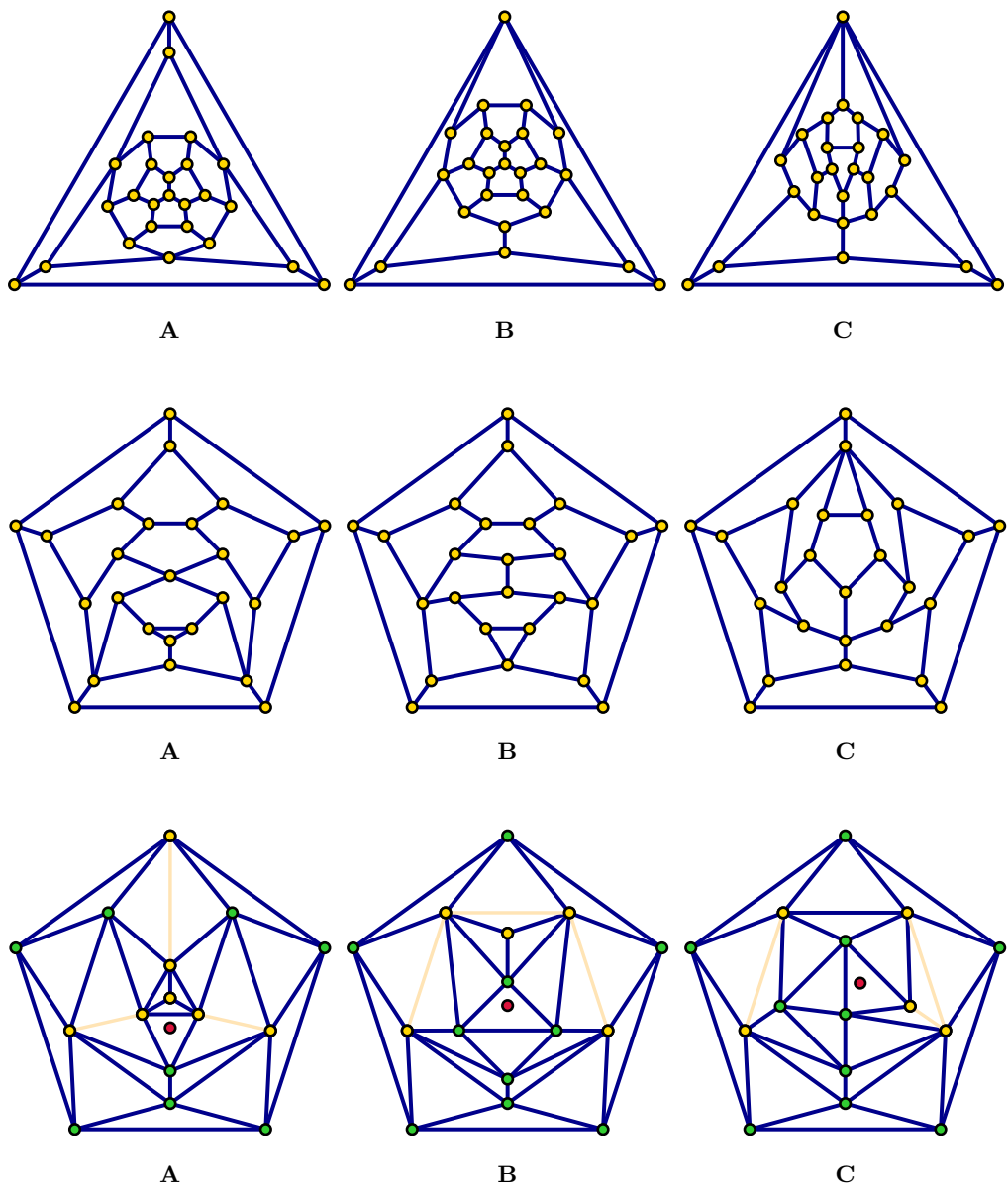

    \centering
    \TriangleMPA{24mm}{0.7mm}{\footnotesize A}
    \TriangleMPB{24mm}{0.7mm}{\footnotesize B}
    \TriangleMPC{24mm}{0.7mm}{\footnotesize C}
    
    \vspace{6mm}

    \TriangleNewMPA{21.8mm}{0.7mm}{\footnotesize A}
    \TriangleNewMPB{21.8mm}{0.7mm}{\footnotesize B}
    \TriangleNewMPC{21.8mm}{0.7mm}{\footnotesize C}

    \vspace{6mm}

    \TriangleTA{21.8mm}{0.7mm}{\footnotesize A}
    \TriangleTB{21.8mm}{0.7mm}{\footnotesize B}
    \TriangleTC{21.8mm}{0.7mm}{\footnotesize C}
    \caption[Pentagulations of a triangle]
    {Three minimal  $3$-connected pentagulations of a triangle\\
    and the corresponding plane triangulations.}
    \label{fig:Triangle}
\end{figure}

Up to isomorphism, there are exactly three $3$-connected pentagulations for $\Delta^+=3$ and $n=3$.
These pentagulations are depicted in Figure \ref{fig:Triangle}.
Two different representations are provided for each of the three pentagulations of the triangle.
In one representation, the leading vertex of the triangulation has degree $3$; 
in the other, the vertex of degree $5$ is chosen as the leading vertex.
The associated plane triangulations are likewise shown in Figure \ref{fig:Triangle}..
For the notation, refer to Remark \ref{remark:colors and lieders}.
It is worth noting that each pentagulation generally corresponds to more than one plane triangulation.
Here we give those that appear more or less symmetrical.
The fact that a pentagulation corresponds to several plane triangulations does not contradict Whitney's theorem. The reason is that the dual of a pentagulation is not a plane triangulation. By adding a few edges (three in this case), we turn this planar graph into a plane triangulation. However, the way these edges are added is not unique, which is why there can be multiple triangulations.
Finally, note that the three pentagulations (A, B, and C) of a triangle are pairwise nonisomorphic:
in pentagulation A, all vertices of the triangular face have degree $3$;
in B, one vertex has degree $4$;
and in C, one vertex has degree $5$.

Up to isomorphism, for $\Delta^+=3$ and $n=7$ there are also exactly three $3$-connected pentagulations.
These pentagulations together with their corresponding plane triangulations are shown in Figure \ref{fig:Heptagon}. 
The pentagulations A, B, and C of a heptagon are pairwise non-isomorphic:
The pentagulation B contains a vertex of degree $5$, unlike A and C, which have no such vertices.
In addition, in C, two vertices of the heptagonal face are adjacent to vertices of degree $4$, 
whereas in A, only one such vertex is adjacent to the heptagonal face.

\begin{figure}[ht]
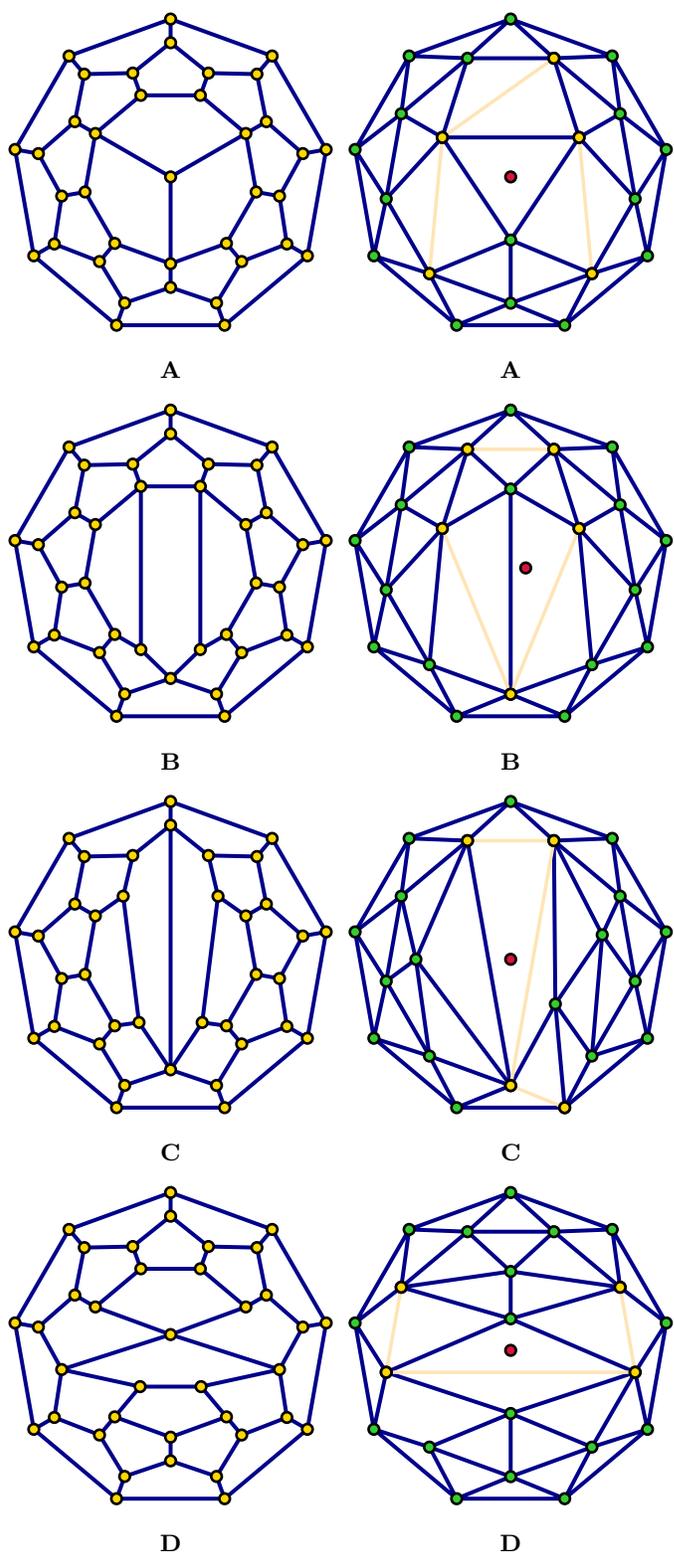

    \centering
    \NonagonMPA{21mm}{0.7mm}{\footnotesize A}
    \NonagonTA{21mm}{0.7mm}{\footnotesize A}

    \NonagonMPB{21mm}{0.7mm}{\footnotesize B}
    \NonagonTB{21mm}{0.7mm}{\footnotesize B}

    \NonagonMPC{21mm}{0.7mm}{\footnotesize C}
    \NonagonTC{21mm}{0.7mm}{\footnotesize C}

    \NonagonMPD{21mm}{0.7mm}{\footnotesize D}
    \NonagonTD{21mm}{0.7mm}{\footnotesize D}

    \caption[Pentagulations of a nonagon]
    {Four minimal  $3$-connected pentagulations of a nonagon\\
    and the corresponding plane triangulations.}
    \label{fig:Nonagon}
\end{figure}

For $\Delta^+=3$ and $n=9$ there are exactly four $3$-connected pentagulations, up to isomorphism.
Figure \ref{fig:Nonagon} depicts these four pentagulations 
alongside their corresponding plane triangulations.
The verification that pentagulations A, B, C, and D are pairwise non-isomorphic is left to the reader 
as a straightforward exercise.

For $n=10$, there are exactly three $3$-connected pentagulations with $\Delta^+=3$, up to isomorphism.
These pentagulations, along with their corresponding plane triangulations, are shown in Figure \ref{fig:Decagon}.
The pentagulations A, B, and C of the decagon are pairwise non-isomorphic,
because the minimum distances from the vertex of degree $6$
to the vertices of the decagonal face are $2$, $1$, and $0$, respectively.

\begin{figure}[ht]
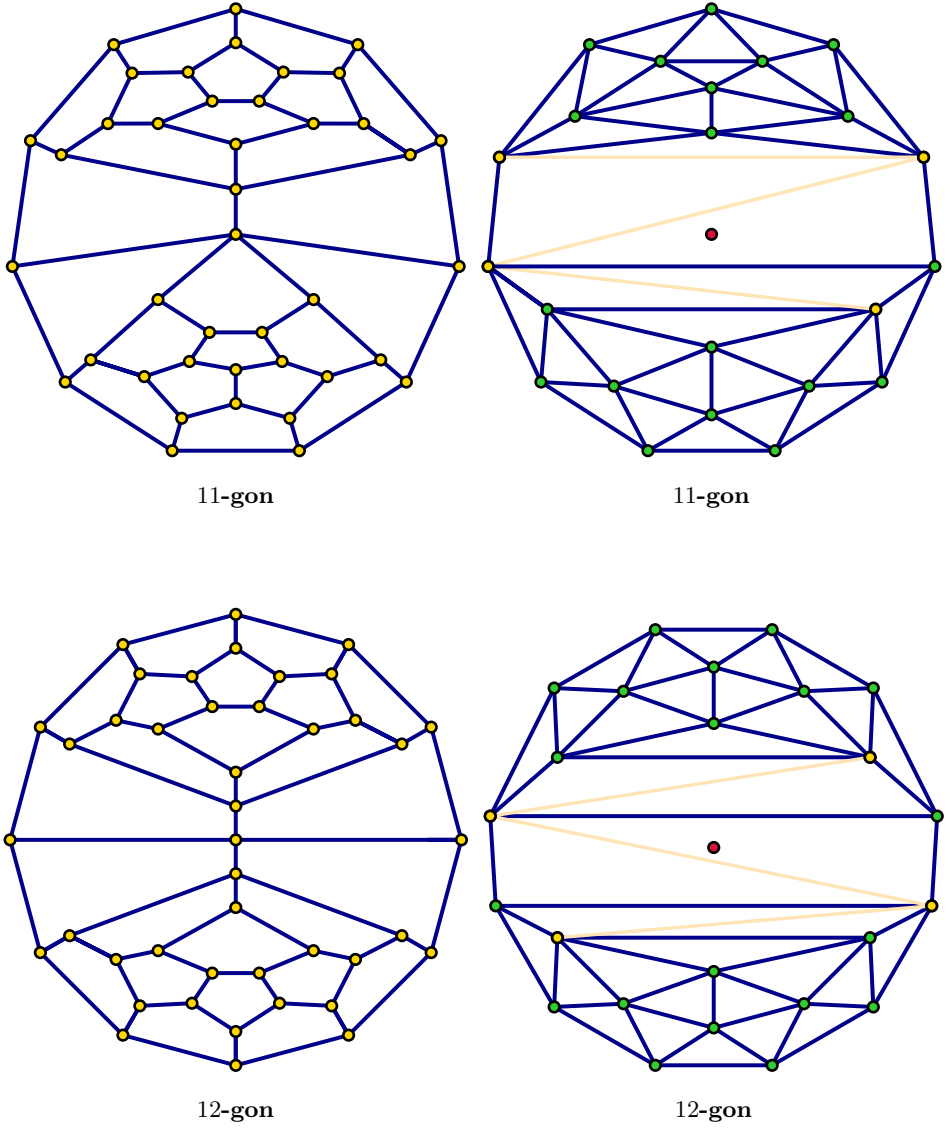

    \centering
    \UndecagonMP{30mm}{0.7mm}{\footnotesize $11$-gon}
    \UndecagonT{30mm}{0.7mm}{\footnotesize $11$-gon}

\vspace{1cm}

    \DodecagonMP{30mm}{0.7mm}{\footnotesize $12$-gon}
    \DodecagonT{30mm}{0.7mm}{\footnotesize $12$-gon}
    \caption[Pentagulations of a undecagon and of a dodecagon]
    {Minimal  $3$-connected pentagulations of an undecagon and\\ a dodecagon
    and the corresponding plane triangulations.}
    \label{fig:Undecagon}
\end{figure}

There is a unique $3$-connected pentagulation (up to isomorphism) for $n=11$ and $n=12$ with $\Delta^+=3$.
These unique pentagulations and their corresponding plane triangulations are depicted 
in Figure \ref{fig:Undecagon}.

Theorem is proved.

We conclude this section by noting that Theorems \ref{Th:exact_scores} and \ref{Th:single} from the introduction
are a consequence of Theorems \ref{Th:dfc>1},
\ref{Th:Triangulations with dfc=2} and \ref{Th:Triangulations with dfc=3}.

\section{An upper bound for \texorpdfstring{$Pg(n)$}{Pg(n)}}
\label{section:Upper bound}
\begin{figure}[ht]
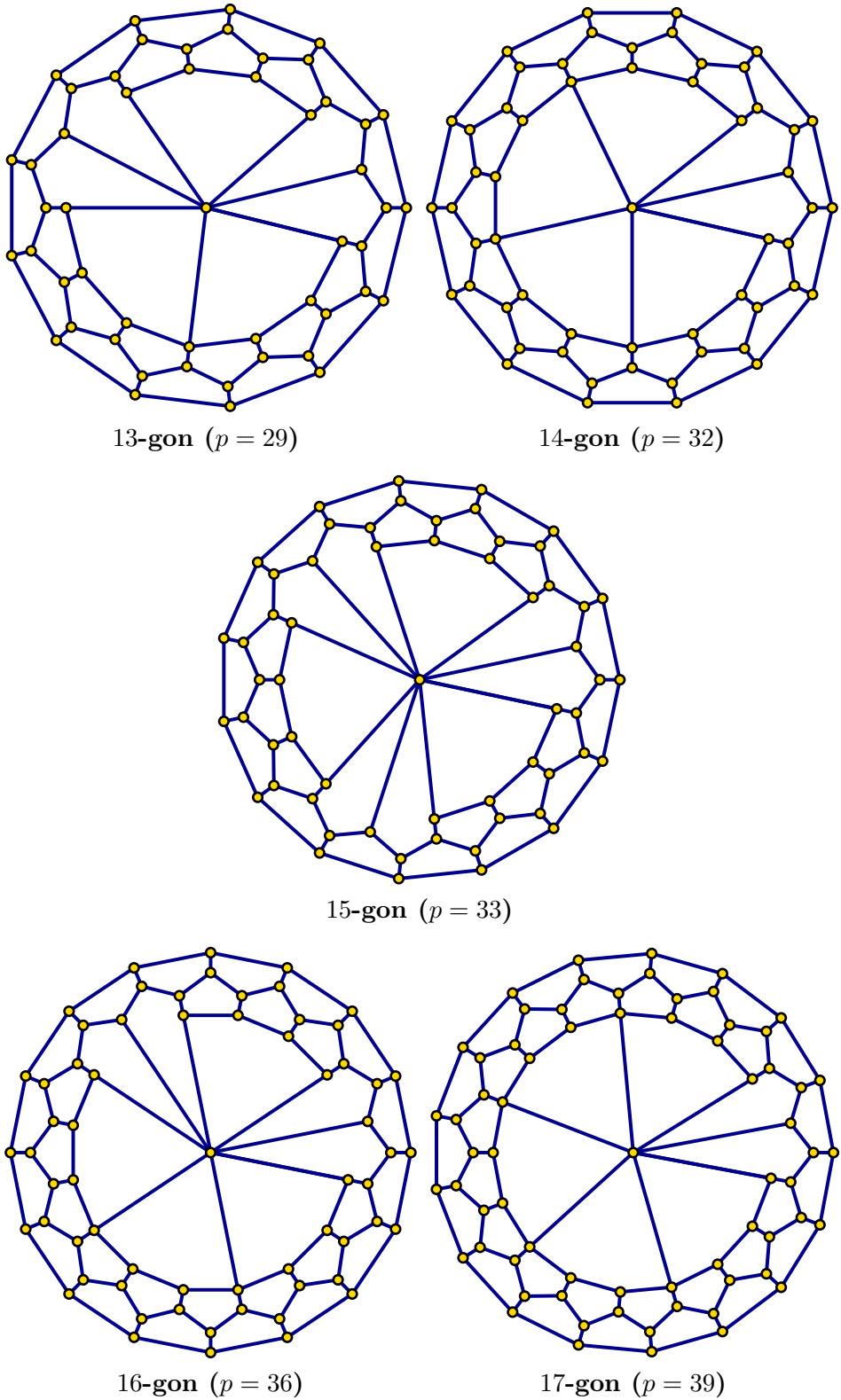

    \centering
    \Bigngons{30mm}{0.7mm}{13}
    \Bigngons{30mm}{0.7mm}{14}

    \Bigngons{30mm}{0.7mm}{15}

    \Bigngons{30mm}{0.7mm}{16}
    \Bigngons{30mm}{0.7mm}{17}
    \caption[Examples of pentagulations for n-gons, $n\in\{13,\ldots,21\}$]
    {Examples of pentagulations for $n$-gons, $n\in\{13,\ldots,17\}$}
    \label{fig:Pentagulations of Big_n-gon}
\end{figure}

In this section, we prove Theorem \ref{Th:upper_bound} stated in the Introduction.
First, note that any integer $n\geq8$ can be written as $n=5k+3l$
with integers $k\geq0$ and $0\leq l\leq4$.
For $n<12$, the representation is as follows:
$$
8=5+3,\ 9=0+3\cdot3,\ 10=2\cdot5+0,\hbox{ and } 11=5+2\cdot3.
$$
For $n\geq12$, one of the five non-negative integers $n,n-3,n-6,n-9,n-12$ must be divisible by $5$.
If $n-3l=5k$, then $n=5k+3l$.
For simplicity, assume $n=5k+3l\geq13$, implying $k\geq1$.
Now we define the graph $G_n$ as follows.
The vertex set of $G_n$ is given by
$$
V(G_n)=\{x_i,y_i,z_i,u_i,v\mid i=0,\ldots,n-1\}\setminus\{u_{5i}\mid i=0,\ldots,k-1\},
$$
and the edges of $G_n$ are defined as follows:
\begin{align*}
    x_ix_{i+1},y_iz_i,z_iy_{i+1},x_iy_i, &&& i=0,\ldots,n-1;\\
    u_iz_i,&&& i=0,\ldots,n-1 \hbox{ and }i\neq5j,j=0,\ldots,k-1;\\
    vz_i,vu_{i\pm1},  &&& i=5j,j=0,\ldots,k-1;\\
    vu_i,    &&& i=5(k-1)+3j+1,j=1,\ldots,l;\\
    u_iu_{i+1},&&& i=0,\ldots,n-1 \hbox{ and }i\neq5j,i\neq5j-1,j=0,\ldots,k-1
\end{align*}
Indices are taken modulo $n$;
for example, $u_n=u_0$, $u_{n+1}=u_1$, etc.
Graphs $G_n$ for specific values of $n$ ($13$, $14$, $15$, $16$, and $17$)
are depicted in Figure \ref{fig:Pentagulations of Big_n-gon}.
The graph $G_n$ includes $n$-gonal face $x_0x_1\ldots x_{n-1}$
and the following pentagonal faces:
\begin{align*}
    z_iy_ix_ix_{i+1}y_{i+1},  &&& i=0,\ldots,n-1;\\
      tu_iz_iy_{i+1}z_{i+1},  &&& i=0,\ldots,n-1,\hbox{ where $t=v$ if $i=5j,j=0,\ldots,k-1$,} \\
                              &&& \hbox{and $t=u_i$ otherwise;}\\
    vu_iu_{i+1}u_{i+2}u_{i+3},&&& i=5j,j=0,\ldots,k-1;\ i=5(k-1)+3j+1,j=0,\ldots,l.
\end{align*}
The graph $G_n$ has
$4n-k+1$ vertices,
$6n+l$ edges, and
$2n+k+l+1$ faces.
Thus, $G_n$ is a pentagulation of an $n$-gon with $2n+k+l$ pentagonal faces, implying $\Pg(n)\leq2n+k+l$.
This completes the proof of Theorem \ref{Th:upper_bound}.

\section{Questions}
\label{section:Questions}
The first and most general open question is:
\begin{question}
How many pentagons are contained in a minimal $3$-connected pentagulation of an $n$-gon  when
$n\geq 13$? In other words, find $\Pg(n)$ for $n\geq13$.
\end{question}
A more specific question is worth highlighting:
\begin{question}
Find $\Pg(13)$.
\end{question}
According to Theorem \ref{Th:exact_scores}, $\Pg(13)$ is at least $27$; on the other hand,
by Theorem \ref{Th:upper_bound}, $\Pg(13)$ is at most $29$. 
Hence, $\Pg(13)$ must be $27$ or $29$,
as the number of pentagons is odd for odd $n$ and even for even $n$.

As we have already mentioned, most of the results were obtained using
computer calculations with the \texttt{plantri} package.
However, there is a slim possibility that these results, 
and possibly stronger ones, could be achieved without a computer.
In this context, a challenge arises:
\begin{question}
Without using a computer, find a proof of these theorems or stronger ones (and possibly even weaker ones).
\end{question}

\newpage
\section{Appendix: Pentagulation diagrams}\label{section:Appendix}

\begin{figure}[ht]
    \centering
    \HeptagonMPA{21.5mm}{0.7mm}{\footnotesize A}
    \HeptagonMPB{21.5mm}{0.7mm}{\footnotesize B}
    \HeptagonMPC{21.5mm}{0.7mm}{\footnotesize C}

\vspace{1cm}

    \HeptagonTA{21.5mm}{0.7mm}{\footnotesize A}
    \HeptagonTB{21.5mm}{0.7mm}{\footnotesize B}
    \HeptagonTC{21.5mm}{0.7mm}{\footnotesize C}
    \caption[Pentagulations of a heptagon]
    {Three minimal  $3$-connected pentagulations of a heptagon\\
    and the corresponding plane triangulations.}
    \label{fig:Heptagon}
\end{figure}

\begin{figure}[ht]
    \centering
    \DecagonMPA{21mm}{0.7mm}{\footnotesize A}
    \DecagonMPB{21mm}{0.7mm}{\footnotesize B}
    \DecagonMPC{21mm}{0.7mm}{\footnotesize C}

\vspace{1cm}

    \DecagonTA{21mm}{0.7mm}{\footnotesize A}
    \DecagonTB{21mm}{0.7mm}{\footnotesize B}
    \DecagonTC{21mm}{0.7mm}{\footnotesize C}
    \caption[Pentagulations of a decagon]
    {Three minimal  $3$-connected pentagulations of a decagon\\
    and the corresponding plane triangulations.}
    \label{fig:Decagon}
\end{figure}

\cleardoublepage

\end{document}